\documentclass[12pt]{article}
\usepackage{amssymb}
\usepackage{a4}
\begin{document}
\def\st{\, : \,}
\def\kbar{{\mathchar'26\mkern-9muk}}  
\def\bra#1{\langle #1 \vert}
\def\ket#1{\vert #1 \rangle}
\def\vev#1{\langle #1 \rangle}
\def\ad{\mbox{ad}\,}
\def\ker{\mbox{Ker}\,}
\def\im{\mbox{Im}\,}
\def\der{\mbox{Der}\,}
\def\ad{\mbox{ad}\,}
\def\b#1{{\mathbb #1}}
\def\c#1{{\cal #1}}
\def\pt{\partial_t}
\def\px{\partial_1}
\def\bpx{\bar\partial_1}
\def\1{{\bf 1}}
\def\la{\langle}
\def\ra{\rangle}
\def\nn{\nonumber \\}
\def\pn{\par\noindent}
\def\etal{{\it et al.}\ }
\def\sq{\mbox{\rlap{$\sqcap$}$\sqcup$}}

\def\R{{\cal R}\,}
\newcommand{\tr}{\triangleright\,}
\newcommand{\tl}{\,\triangleleft}
\newcommand{\tro}{\triangleright^{op}\,}
\newcommand{\tlo}{\,^{op}\triangleleft}
\def\cross{{\triangleright\!\!\!<}}
\def\cocross{{>\!\!\!\triangleleft}}
\def\g{\mbox{\bf g\,}}
\def\uqg{\mbox{$U_q{\/\mbox{\bf g}}$ }}
\def\uqs{\mbox{$U_q{\/\mbox{so}(N)}$}}
\def\uqgo{\mbox{$U^{op}_q{\/\mbox{\bf g}}$ }}
\def\uqgp{\mbox{$U_q^+{\/\mbox{\bf g}}$ }}
\def\uqgn{\mbox{$U_q^-{\/\mbox{\bf g}}$ }}

\renewcommand{\thefootnote}{\fnsymbol{footnote}}

\newcommand{\sect}[1]{\setcounter{equation}{0}\section{#1}}
\renewcommand{\theequation}{\thesection.\arabic{equation}}
\newcommand{\subsect}[1]{\setcounter{equation}{0}\subsection{#1}}
\newcommand{\app}[1]{\setcounter{section}{0}
\setcounter{equation}{0} \renewcommand{\thesection}
{\Alph{section}}\section{#1}}

\newcommand{\be}{\begin{equation}}
\newcommand{\ee}{\end{equation}}
\newcommand{\bea}{\begin{eqnarray}}
\newcommand{\eea}{\end{eqnarray}}
\newcommand{\ba}{\begin{array}}
\newcommand{\ea}{\end{array}}
%
%
%
\newtheorem{prop}{Proposition}
\newtheorem{lemma}{Lemma}
\newtheorem{theorem}{Theorem}
\newtheorem{corollary}{Corollary}
%
%
%
\newenvironment{proof}[1]{\vspace{5pt}\noindent{\bf Proof #1}\hspace{6pt}}%
{\hfill\sq}
\newcommand{\bp}{\begin{proof}}
\newcommand{\ep}{\end{proof}\par\vspace{10pt}\noindent}
%
%

\title{Unbraiding the braided tensor product}

\author{Gaetano Fiore,$\strut^{1,2}$ \, Harold Steinacker$\strut^{3}$,
        Julius Wess$\strut^{3,4}$ \\\\
        \and
        $\strut^1$Dip. di Matematica e Applicazioni, Fac.  di Ingegneria\\ 
        Universit\`a di Napoli, V. Claudio 21, 80125 Napoli
        \and
        $\strut^2$I.N.F.N., Sezione di Napoli,\\
        Complesso MSA, V. Cintia, 80126 Napoli
        \and
        $\strut^3$Sektion Physik, Ludwig-Maximilian Universit\"at,\\
        Theresienstra\ss e 37, D-80333 M\"unchen
        \and
        $\strut^4$Max-Planck-Institut f\"ur Physik\\
        F\"ohringer Ring 6, D-80805 M\"unchen
        }
\date{}

\maketitle
\abstract{We show that the braided tensor product algebra 
$\c{A}_1\underline{\otimes}\c{A}_2$ of two module
algebras $\c{A}_1, \c{A}_2$ of a quasitriangular Hopf algebra $H$ 
is equal to the ordinary tensor product algebra of $\c{A}_1$ with
a subalgebra of $\c{A}_1\underline{\otimes}\c{A}_2$
isomorphic to $\c{A}_2$, provided there exists a 
realization of $H$ within $\c{A}_1$. In other words, under this 
assumption we construct a transformation of generators which `decouples' 
$\c{A}_1, \c{A}_2$ (i.e. makes them commuting).
We apply the theorem to the braided tensor product 
algebras of two or more quantum group covariant quantum spaces, deformed 
Heisenberg algebras and $q$-deformed fuzzy spheres.}

\newpage

\sect{Introduction and main theorem}

As is well known, given two associative unital algebras 
$\c{A}_1, \c{A}_2$ (over the field $\b{C}$, say), 
one can build a new module algebra $\c{A}$ which is 
as a vector space the tensor product $\c{A}=\c{A}_1\otimes\c{A}_2$
of the two vector spaces (over the same field)  by postulating the 
product law
\be
(a_1\otimes a_2)(b_1\otimes b_2)=a_1b_1\otimes a_2b_2. \label{primo}
\ee
The resulting algebra is the ordinary tensor product algebra.
(\ref{primo}) is equivalent to the set of relations
\bea
(a_1\otimes \1_2)(b_1\otimes \1_2)&=& a_1b_1\otimes \1_2,\label{trivial1} \\
(a_1\otimes \1_2)(\1_1\otimes a_2)&=& a_1\otimes a_2,\label{trivial2} \\
(\1_1\otimes a_2)(\1_1\otimes b_2)&=& \1_1\otimes a_2b_2, \label{trivial3}\\
(\1_1\otimes a_2)(a_1\otimes \1_2)&=&(a_1\otimes \1_2)(\1_1\otimes a_2). 
\label{trivial4}
\eea
However, in many cases
the same goal can be reached also by replacing (\ref{trivial4}) 
by some suitable nontrivial commutation relations.
With a standard abuse of notation we shall denote in the sequel
$a_1\otimes a_2$ by $a_1a_2$ for any $a_1\in\c{A}_1$, 
$a_2\in\c{A}_2$ and omit all units $\1_i$ when multiplied
by non-unit elements; consequently (\ref{trivial1}-\ref{trivial3})
take trivial forms, whereas  (\ref{trivial4}) becomes the
commutation relation
\be
a_2a_1=a_1a_2.                                       \label{trivial}
\ee

If $\c{A}_1, \c{A}_2$ are module algebras of a Lie algebra
$\g$, and we require $\c{A}$ to be too, then (\ref{trivial})
has no alternative, because any 
$g\in\g$ acts as a derivation on the (algebra as well as tensor)
product of any two elements or,  in Hopf algebra language,
because the coproduct  $\Delta(g)=g_{(1)}\otimes g_{(2)}$ 
(at the rhs we have used Sweedler notation)
of the Hopf algebra $H\equiv U\g$ is cocommutative.
In the main part of this paper we shall work with 
right-module algebras (instead of left ones), and denote by 
$\tl:(a_i,g)\in\c{A}_i\times H\rightarrow a_i\tl g\in\c{A}_i$ 
the right action; the reason is that they are equivalent 
to left comodule algebras, which 
are used in much of the literature. In section \ref{leftsec}
we shall give the formulae for left module algebras.
We recall that a right action 
$\tl:(a,g)\in\c{A}\times H\rightarrow a\tl g\in\c{A}$ 
by definition fulfills
\begin{eqnarray}
&&a\tl(gg') = (a\tl g) \tl   g',                     \label{modalg1}\\
&&(aa')\tl g = (a\tl g_{(1)})\, (a'\tl g_{(2)}).       \label{modalg2}
\end{eqnarray}

If we ``$q$-deform'' this setting by taking as Hopf algebra
$H$ the quantum group $\uqg$, and as
$\c{A}_i$ the corresponding $q$-deformed module algebras, then
it is also known \cite{JoyStr86,Maj1,Maj95} that
although $\Delta(g)$ is no longer cocommutative,
it is still possible to build
the deformed counterpart of $\c{A}$ if one replaces  
(\ref{trivial}) with nontrivial
commutation relations of the form
\be
a_2a_1= (a_1\tl \R^{(1)})\, (a_2\tl \R^{(2)}).     \label{absbraiding}
\ee
Here $\R\equiv \R^{(1)}\otimes \R^{(2)}\in H^+\otimes H^-$ 
denotes the socalled universal $R$-matrix of $H\equiv\uqg$ \cite{Dri86},
and $H^{\pm}$ denote the Hopf positive and negative 
Borel subalgebras of $H$.
This yields instead of $\c{A}$ a {\it braided} tensor product algebra
$\c{A}^+=\c{A}_1\underline{\otimes}^+\c{A}_2$ \cite{Maj95}.
An alternative one $\c{A}^-=\c{A}_1\underline{\otimes}^-\c{A}_2$
is obtained by
replacing in the previous formula $\R$ by $\R^{-1}_{21}$: 
\be
a_2a_1= (a_1\tl \R^{-1(2)})\, (a_2\tl \R^{-1(1)}).   \label{absbraiding'}
\ee
Both $\c{A}^+$ and $\c{A}^-$ go to the ordinary tensor product algebra
$\c{A}$ in the limit $q\to 1$.

This is a particular example of a more general notion, that of a
{\it crossed (or twisted) tensor product} \cite{BorMar00}
of two unital associative algebras.

In view of (\ref{absbraiding}) or (\ref{absbraiding'}) 
studying representations of $\c{A}^{\pm}$ is a more difficult
task than just studying the representations of $\c{A}_1,\c{A}_2$,
taking their tensor products and studying
the irreducible ones there contained. The degrees of freedom
of $\c{A}_1,\c{A}_2$ are so to say ``coupled''. One might ask 
whether one can ``decouple'' them by a transformation of generators.

In this work we present a sufficient condition 
for the construction of a transformation making
$\c{A}^+$ equal to an ordinary
tensor product $\c{A}_1\otimes\tilde\c{A}_2^+$, 
with $\tilde\c{A}_2^+$ a subalgebra 
of $\c{A}^+$ isomorphic to $\c{A}_2$ and {\it commuting} 
[in the sense (\ref{trivial})] with $\c{A}_1$, 
although - of course - no longer a $H$-submodule; and
similarly for $\c{A}^-$.  In a quantum theory framework 
one could thus interpret the generators of 
$\c{A}_1,\tilde\c{A}_2^{\pm}$ as
pertaining to decoupled degrees of freedom, describing
e.g. some composite or ``quasiparticle'' excitations. 
Reducing $\c{A}^{\pm}$ to a form
$\c{A}_1\otimes\tilde\c{A}_2^{\pm}$ will be called an
{\it unbraiding} of the braided tensor product algebra
$\c{A}^{\pm}=\c{A}_1\underline{\otimes}^{\pm}\c{A}_2$.
The sufficient condition is that
there respectively exists an algebra homomorphism $\varphi_1^+$ or
an algebra homomorphism $\varphi_1^-$
\be
\varphi_1^{\pm}: \c{A}_1\cocross H^{\pm} \rightarrow \c{A}_1
                                            \label{Hom}
\ee
acting as the identity on $\c{A}_1$, namely for any $a_1\in \c{A}_1$
\be
\varphi^{\pm}_1(a_1)=a_1  .                          \label{ident0}
\ee
(Note that, as a consequence of (\ref{ident0}), 
$\varphi^{\pm}_1$ is idempotent,
$(\varphi^{\pm}_1)^2=\varphi^{\pm}_1$).
Here $\c{A}_1\cocross H^{\pm}$ denotes the cross product between 
$\c{A}_1$ and $H^{\pm}$. In other words, this amounts
to assuming that $\varphi_1^+(H^+)$ [resp. $\varphi_1^-(H^-)$] provides a
{\it realization} of $H^+$ (resp. $H^-$) within $\c{A}_1$.
In fact, $\tilde\c{A}_2$ is found using the main result of this
work: 
\begin{theorem} 
Let $\{H,\R\}$ be a quasitriangular Hopf algebra and 
$H^+,H^-$ be Hopf subalgebras of $H$ such that $\R\in H^+\otimes H^-$.
Let $\c{A}_1, \c{A}_2$ be respectively a $H^+$- and
a $H^-$-module algebra, so that we can define $\c{A}^+$ 
as in (\ref{absbraiding}), and $\varphi_1^+$ be a homomorphism
of the type (\ref{Hom}), (\ref{ident0}), so that we can define
the ``unbraiding'' 
map $\chi^+:\c{A}_2\rightarrow \c{A}^+$ by
\be
\chi^+(a_2):= \varphi_1^+(\R^{(1)})\, (a_2\tl \R^{(2)}). 
\label{Def+}
\ee
Alternatively, let $\c{A}_1, \c{A}_2$ be respectively a $H^-$- and
a $H^+$-module algebra, so that we can define $\c{A}^-$ 
as in (\ref{absbraiding'}), and $\varphi_1^-$ be a homomorphism
of the type (\ref{Hom}), (\ref{ident0}), so that we can define
the ``unbraiding'' 
map $\chi^-:\c{A}_2\rightarrow \c{A}^-$ by
\be
\chi^-(a_2):= \varphi_1^-(\R^{-1(2)})\, (a_2\tl \R^{-1(1)}). 
\label{Def-}
\ee
In either case $\chi^{\pm}$ are then injective algebra homomorphisms and 
\be
[\chi^{\pm}(a_2),\c{A}_1]=0,                       \label{commu}
\ee
namely the subalgebras 
$\tilde\c{A}_2^{\pm}:=\chi^{\pm}(\c{A}_2)\approx\c{A}_2$ commute
with $\c{A}_1$. Moreover $\c{A}^{\pm}=\c{A}_1\otimes\tilde\c{A}_2^{\pm}$.
\label{maintheo}
\end{theorem}
\bp{}. We start by recalling the content of 
the hypotheses stated in the theorem.
The algebra $\c{A}_1\cocross H^{\pm}$ as a vector space is the tensor
product of $\c{A}_1$ and $H^{\pm}$, whereas its product law is
obtained combining the product laws of these two tensor factors with the
cross-product law,
\be
a_1 g=g_{(1)}\, (a_1\tl g_{(2)}),                          \label{crossprod}
\ee
for any $a_1\in\c{A}_1$ and $g\in H^{\pm}$.
$\varphi_1^{\pm}$ being an algebra homomorphism means that for any 
$\xi,\xi'\in\c{A}_1\cocross H^{\pm}$
\be
\varphi^{\pm}_1(\xi\xi')=\varphi^{\pm}_1(\xi)\,\varphi^{\pm}_1(\xi').      
   \label{HOM}
\ee
For $\xi\equiv a\in\c{A}_1\subset\c{A}_1\cocross H^{\pm}$, 
$\xi'\equiv g\in H^{\pm}\subset\c{A}_1\cocross H^{\pm}$ 
this implies
\be
a\varphi^{\pm}(g)=\varphi^{\pm}(g_{(1)}) (a\tl g_{(2)})        \label{Homo1}
\ee
Hereby we have also used (\ref{ident0}) and (\ref{crossprod}).
After these preliminaries, note that 
under the assumption (\ref{absbraiding}),
for any $a_1\in\c{A}_1$ and  $a_2\in\c{A}_2$
\bea
a_1\,\chi^+(a_2) & \stackrel{(\ref{Def+})}{=}& a_1
\varphi_1^+(\R^{(1)})\, (a_2\tl \R^{(2)}) \nn
& \stackrel{(\ref{Homo1})}{=}& \varphi^+(\R^{(1)}_{(1)}) 
(a_1\tl \R^{(1)}_{(2)}) (a_2\tl \R^{(2)})\nn
& \stackrel{(\ref{delta1})}{=}& \varphi^+(\R^{(1)}) 
(a_1\tl \R^{(1')}) (a_2\tl \R^{(2)}\R^{(2')})\nn
& \stackrel{(\ref{absbraiding})}{=}& \varphi^+(\R^{(1)}) 
(a_2\tl \R^{(2)}) a_1\nn
& \stackrel{(\ref{Def+})}{=}&\chi^+(a_2) \, a_1, \nonumber
\eea
which proves (\ref{commu}) in this case. Moreover
\bea
\chi^+(a_2a_2') & \stackrel{(\ref{Def+})}{=}&
\varphi_1^+(\R^{(1)})\, (a_2a_2'\tl \R^{(2)}) \nn
& \stackrel{(\ref{modalg2})}{=}& \varphi^+(\R^{(1)}) 
(a_2\tl \R^{(2)}_{(1)}) (a_2'\tl \R^{(2)}_{(2)})\nn
& \stackrel{(\ref{delta2})}{=}& \varphi^+(\R^{(1')}\R^{(1)}) 
(a_2\tl \R^{(2)}) (a_2'\tl \R^{(2')})\nn
& \stackrel{(\ref{HOM})}{=}& \varphi^+(\R^{(1')}) \varphi(\R^{(1)}) 
(a_2\tl \R^{(2)}) (a_2'\tl \R^{(2')})\nn
& \stackrel{(\ref{Def+})}{=}& \varphi^+(\R^{(1')}) \chi^+
(a_2) (a_2'\tl \R^{(2')})\nn
& \stackrel{(\ref{commu})}{=}&  \chi^+
(a_2) \varphi(\R^{(1')})(a_2'\tl \R^{(2')})\nn
& \stackrel{(\ref{Def+})}{=}&\chi^+(a_2) \, \chi^+(a_2'), \nonumber
\eea
proving that $\chi^+$ is a homomorphism. To prove injectivity we 
show that $\chi^+$ can be inverted on $\chi^+(\c{A}_2)$, and
the inverse is given by
\be
(\chi^+)^{-1}(\tilde a_2)=V^{-1}\,
\left([\varphi^+(S^{-1}\R^{(1)})\tilde a_2]\tl \R^{(2)}\right)
                                                        \label{inve}
\ee
where $V\in\c{A}_1$ is the invertible element
defined by $V:= \varphi_1^+(S^{-1}\R^{(1)})\tl \R^{(2)}$
($V$ is invertible because $\R$ is). In fact,
\bea
 & & V^{-1}\,[\varphi_1^+(S^{-1}\R^{(1)})\chi^+(a_2)]\tl
 \R^{(2)} \nn
&\stackrel{(\ref{Def+})}{=} & V^{-1}\,
[\varphi^+_1(S^{-1}\R^{(1)})
\varphi_1^+(\R^{(1')})\, (a_2\tl \R^{(2')})]\tl \R^{(2)} \nn
& \stackrel{(\ref{SR}),(\ref{HOM})}{=}& V^{-1}\, \left\{\varphi^+_1
[S^{-1}(\R^{-1(1')}\R^{(1)})]\, (a_2\tl \R^{-1(2')})
\right\}\tl \R^{(2)} \nn
& \stackrel{(\ref{modalg2})}{=}& 
V^{-1}\, \varphi^+_1[S^{-1}(\R^{-1(1')}\R^{(1)})
]\tl \R^{(2)}_{(1)}\:\: (a_2\tl \R^{-1(2')})\tl \R^{(2)}_{(2)} \nn
& \stackrel{(\ref{delta2})}{=}& 
V^{-1}\, \varphi^+_1[S^{-1}(\R^{-1(1')}\R^{(1)}\R^{(1")})
]\!\tl\! \R^{(2")}\, (a_2\!\tl \R^{-1(2')})\!\tl\! \R^{(2)}) \nn
& \stackrel{(\ref{modalg1})}{=}& 
V^{-1}\, \varphi^+_1[S^{-1}(\R^{(1")})
]\tl \R^{(2")}\, a_2\nn
&=& V^{-1}\, V a_2=a_2.         \nonumber
\eea
In fact,
if $\varphi_1^+$ can be extended to an algebra homomorphism
$\varphi_1: \c{A}_1\cocross H \rightarrow \c{A}_1$ a little
calculation with the help of eq.'s (\ref{int}), (\ref{defu}) shows that
$V=\varphi_1(v)$, where $v\in H$ is the invertible central element
defined by (\ref{Defv}).  We know that 
$\c{A}_1\otimes\tilde\c{A}_2^+\subset \c{A}^+$.
To prove that $\c{A}^+=\c{A}_1\otimes\tilde\c{A}_2^+$ note first that 
by (\ref{absbraiding}) any element in $\c{A}^+$ can be written 
as a sum of products $a_1a_2$, with $a_1\in \c{A}_1$ and 
$a_2\in\c{A}_2$. So we need to show that
\be
a_1a_2=b^{(1)}\chi^+(b^{(2)})                       \label{ciccio}
\ee
for some $b^{(1)}\in \c{A}_1$, $b^{(2)}\in \c{A}_2$ (at the rhs
a sum of many terms is implicitly understood). Now this can be
proved as follows:
\bea
a_1a_2 & = & a_1\varphi_1^+(\1_H)\:
(a_2\tl \1_H)= a_1\varphi_1^+(\R^{-1(1)}\R^{(1')})
\:[a_2\tl  (\R^{-1(2)} \R^{(2')})]\nn
&\stackrel{(\ref{modalg1})}{=} & a_1\varphi_1^+(\R^{-1(1)})
\varphi_1^+(\R^{(1')})\:
(a_2 \tl \R^{-1(2)})\tl \R^{(2')}\nn
& \stackrel{(\ref{Def+})}{=} & a_1\varphi_1^+(\R^{-1(1)})\:
\chi^+(a_2 \tl \R^{-1(2)}), \nonumber
\eea
which is of the form (\ref{ciccio}).

The proof for $\chi^-$
under the corresponding assumptions is completely analogous.
\ep
In the next Section we shall need an alternative expression for $\chi^{\pm}$, 
which we prove in the appendix:
\begin{prop}
\bea
&&\chi^+(a_2)=(a_2\tl \R^{-1(2)})\, \varphi^+_1(S\R^{-1(1)} ), 
\label{past1} \\
&&\chi^-(a_2)=(a_2\tl \R^{(1)})\, \varphi^-_1(S\R^{(2)} ).\label{past2} 
\eea
\label{propalt}
\end{prop}

The rest of the paper is essentially devoted to
illustrate the application of
Theorem \ref{maintheo} to some algebras $\c{A}_i$ for which homomorphisms
$\varphi_1^{\pm}$ are known.
In Ref. \cite{CerFioMad00} algebra homomorphisms
$\varphi^{\pm}_1$ have been found
for (a slightly enlarged version $\c{A}_1$ of)
the algebra of functions on the $N$-dimensional quantum 
Euclidean space \cite{FadResTak89} $\b{R}_q^N$, corresponding to $H=U_qso(N)$. 
The explicit forms of $\varphi^{\pm }_1$ on the Faddeev-Reshetikhin-Takhtadjan
(FRT) generators $\c{L}^{\pm}{}^i_j$ of $U_qso(N)$ are recalled in the 
appendix \ref{appeuc}.
The maps $\varphi^{\pm}_1$ for $N=3$ 
are given also in 
Ref. \cite{CerMadSchWes00}. The same maps
do the job also on the quotient spaces obtained by
setting $x^ix_i=1$ [quantum $(N\!-\!1)$-dimensional spheres $S_q^{N-1}$],
and the appropriate maps for the $q$-deformed fuzzy sphere $S_{q,M}^2$ 
have been found in \cite{fuzzyq}. Therefore 
\uqs and the quantum Euclidean spaces/spheres 
provide nontrivial $H$ and $\c{A}_1$ for the application
of the above theorem. In fact, the constructions of the frame  given in
Ref. \cite{FioMad99,CerFioMad00} can be interpreted as an application of the
theorem with $\c{A}_1\equiv \b{R}_q^N$  and $\c{A}_2$
the $N!$-dim exterior algebra generated by the differentials $dx^i$ of
the $\uqs$-covariant differential calculus (although 
with a universal $R$-matrix
$\R$ slightly modified by multiplication by the coproduct 
$\Delta(\Lambda)=\Lambda\otimes\Lambda$ of a new element
$\Lambda$ generating dilatations); 
consequently, in agreement with the
philosophy of Ref. \cite{Mad99c}, 
the algebra of differential forms on $\b{R}_q^N$ can be written as
$\b{R}_q^N\otimes\tilde\c{A}_2$, where $\tilde\c{A}_2$  is
the $N!$-dim exterior algebra generated by the frame elements. 
On the other hand, the existence 
of algebra homomorphisms
$\varphi:\c{A}_1\cocross H\rightarrow \c{A}_1$, for 
$H=\uqs,U_qsl(N)$  and $\c{A}_1$ respectively equal to
(a suitable completion of) the
$\uqs$-covariant Heisenberg algebra or the $U_qsl(N)$-covariant 
Heisenberg or Clifford algebras,
has been known for even a longer time \cite{Fiocmp95,ChuZum95,Hay90}, 
so the theorem also applies if we choose as 
$(H,\c{A}_1)$ one of these pairs of algebras.

Of course the above theorem can be used iteratively to completely
unbraid an algebra $\c{A}$ obtained by repeated braided 
tensor product [through prescription (\ref{absbraiding}), 
or prescription (\ref{absbraiding'})] of an arbitrary 
number of  $H$-module algebras $\c{A}_1,\c{A}_2,...,\c{A}_M$.
We shall explicitly consider the particular case that the 
latter be $M$ identical copies of the $\uqs$-covariant
quantum space/sphere (Section \ref{qspaces}), 
of the $U_qso(3)$-covariant $q$-fuzzy sphere (Section \ref{q_fuzzy}), 
or of the  \uqs- or $U_qsl(N)$-covariant Heisenberg algebra
(Section \ref{heisenberg}). There we shall explicitly write down
the generators of $\tilde\c{A}_2^{\pm}$ for the lowest $N$ examples. 

In appendix \ref{appeuc} we analyze the properties of $\varphi^{\pm}$
under the main real sections of $\uqs$, what was left aside in Ref.
\cite{CerFioMad00}. In section \ref{star} we investigate in the
context of general position the properties of $\chi^{\pm}$ 
under the $*$-structures.

\sect{The unbraiding under the $*$-structures}
\label{star}

Assume $H$ is a Hopf $*$-algebra, namely the coproduct $\Delta$ and
counit $\varepsilon$ are $*$-homomorphisms,
\be
\Delta(g^*)\equiv(g^*)_{(1)}\otimes (g^*)_{(2)}
=(g_{(1)})^*\otimes (g_{(2)})^*,           \label{starhopf}
\ee
and $\c{A}_1$, $\c{A}_2$ are $H$-module
$*$-algebras, namely for any $a_i\in\c{A}_i$
\be
(a_i\tl g)^*=a_i^*\tl S^{-1} g^*                       \label{starmod}
\ee
(here $S$ denotes the antipode of $H$);
we have used and shall use the same symbol $*$ for the $*$-structure
on all algebras $H,\c{A}_1$, etc. Then $*$ is  
a $*$-structure also for $\c{A}_1\cocross H$. 
The same statement is not automatically true for the braided 
tensor product algebra 
$\c{A}^{\pm}=\c{A}_1\underline{\otimes}^{\pm}\c{A}_2$, 
because the basic requirement that the
latter be antimultiplicative
\be
(a_2a_1)^*=a_1^*a_2^*                          \label{antimult}
\ee
(note that this would make $\c{A}^{\pm}$ also
a $H$-module $*$-algebra) is not automatically guaranteed.
In fact, applying this would-be $*$ to
rhs of (\ref{absbraiding}) one finds
\bea
(a_2a_1)^* & \stackrel{(\ref{absbraiding})}{=}& 
\left[(a_1\tl \R^{(1)})\, (a_2\tl \R^{(2)})\right]^*\nn
& \stackrel{(\ref{antimult})}{=}& (a_2\tl \R^{(2)})^* (a_1\tl \R^{(1)})^*\nn
& \stackrel{(\ref{starmod})}{=}& (a_2^*\tl S^{-1}\R^{(2)}{}^*) 
(a_1^*\tl S^{-1}\R^{(1)}{}^*) \nn
& \stackrel{(\ref{absbraiding})}{=}& 
(a_1^*\tl S^{-1}\R^{(1)}{}^*\,\R^{(1')})(a_2^*\tl S^{-1}\R^{(2)}{}^*\,
\R^{(2')});                                     \label{richiamo}
\eea
in order that this be equal to the rhs of (\ref{antimult})
it is necessary that $(S^{-1}\otimes S^{-1})\R^*=\R^{-1}$,
which upon use of (\ref{SR}) is equivalent to
\be
\R^*=\R^{-1}                                       \label{real1}
\ee
(here $\R^* $ means $\R^{(1)}{}^*\otimes \R^{(2)}{}^*$).
This condition is fulfilled only
for the standard noncompact sections (\ref{|q|=1}) of $U_q \g$,
for $|q|=1$; as a consequence, 
$\c{A}^+=\c{A}_1\underline{\otimes}^+\c{A}_2$ becomes a $H$-module $*$-algebra
if one extends the $*$-structures of the tensor factors
to $\c{A}$ using (\ref{antimult}). The same holds for $\c{A}^-$.

On the contrary, the compact section, which requires
$q\in\b{R}$, is characterized by
\be
\R^*=\R_{21}.                                     \label{real2}
\ee
In the latter case the map $*$ introduced 
through (\ref{antimult}) makes sense only as
an involutive antimultiplicative antilinear map
$\c{A}^+\to\c{A}^-$, if both $\c{A}^+$ and $\c{A}^-$ exist.
In fact, in this case the last line in (\ref{richiamo})
will be replaced by
$$
\stackrel{(\ref{absbraiding'}),(\ref{real2})}{=}
(a_1^*\tl S^{-1}\R^{(2)}\,\R^{-1(2')})(a_2^*\tl S^{-1}\R^{(1)}\,
\R^{-1(1')})\stackrel{(\ref{SR})}{=} a_1^*a_2^*,
$$
as required.
Alternatively, if $\c{A}_1,\c{A}_2$ are two copies of the same
algebra and we denote by $\psi:\c{A}_1\to\c{A}_2$ the map
associating to each $a_1\in\c{A}_1$ the equivalent element
in $\c{A}_2$, one can define an alternative $*$-structure
$\star$ in $\c{A}^{\pm}$ by setting
\be
a_1^{\star}=\psi(a_1^*)\qquad a_2^{\star}=\psi^{-1}(a_2^*), 
                                 \label{altstar}
\ee
since this is instead
compatible with (\ref{absbraiding}). In fact, (\ref{richiamo}) 
will become
\bea
(a_2a_1)^{\star} & \stackrel{(\ref{absbraiding})}{=}& 
\left[(a_1\tl \R^{(1)})\, (a_2\tl \R^{(2)})\right]^{\star}\nn
& \stackrel{(\ref{antimult})}{=}& (a_2\tl \R^{(2)})^{\star} 
(a_1\tl \R^{(1)})^{\star}\nn
& \stackrel{(\ref{starmod})}{=}& \psi^{-1}(a_2^*\tl S^{-1}\R^{(2)}{}^*) 
\psi(a_1^*\tl S^{-1}\R^{(1)}{}^*) \nn
& \stackrel{(\ref{absbraiding})}{=}& 
\psi(a_1^*\tl S^{-1}\R^{(2)}\,\R^{-1(2')})\,\psi^{-1}
(a_2^*\tl S^{-1}\R^{(1)}\,\R^{-1(1')})\nn
& \stackrel{(\ref{SR})}{=} &  
\psi(a_1^*)\,\psi^{-1}(a_2^*)a_1\nn
& \stackrel{(\ref{altstar})}{=} &  a_1^{\star}\, a_2^{\star}
\eea
A similar trick can be used also if one considers an
iterated braided tensor product of $M>2$
copies of the same algebra, see next section. However, 
such $\star$'s have not the standard commutative limit,
because of the presence of the map $\psi$.

Inspired by the applications of
the next two Sections, we now assume that $\varphi^{\pm}_1$
fulfill some specific conditions relating its action before
and after the application of the involution $*$, 
and analyze the identities relating the action 
of $\chi^{\pm}$ before
and after the application of $*$  which follow herefrom.

\begin{prop} Assume that the conditions of Theorem
\ref{maintheo} for defining $\chi^+$ or $\chi^-$ are
fulfilled. If $\R^*=\R^{-1}$ and for any $g^{\pm}\in H^{\pm}$
\be
[\varphi_1^{\pm}(g^{\pm})]^*=\varphi_1^{\pm} (g^{\pm}{}^*),            
\label{hypo1}
\ee
in other words $\varphi_1^{\pm}$ are $*$-homomorphisms, then
\be
[\chi^{\pm}(a_2) ]^* =\chi^{\pm}(a_2^*).             \label{result1}
\ee
If $\R^*=\R_{21}$ and $*:H^{\pm}\to H^{\mp}$ fulfills
\be
[\varphi_1^{\pm}(g)]^*=\varphi_1^{\mp} (g^*),            \label{hypo2}
\ee 
then
\be
[\chi^{\pm}(a_2)]^*=\chi^{\mp}(a_2^*).                   \label{result2}
\ee
\label{propstar}
\end{prop}
\bp{}. Under the first assumptions, for any $a_2\in \c{A}_2$, 
\bea
[\chi^+(a_2) ]^* & \stackrel{(\ref{past1})}{=} & 
[(a_2\tl \R^{-1(2)})\, \varphi^+_1(S\R^{-1(1)} )]^* \nn
 &=& [\varphi^+_1(S\R^{-1(1)} )]^* (a_2\tl \R^{-1(2)})^* \nn
& \stackrel{(\ref{hypo1}),(\ref{starmod})}{=} & 
\varphi^+_1(S^{-1}\R^{-1(1)}{}^* ) 
(a_2^*\tl S^{-1}\R^{-1(2)}{}^*) \nn
& \stackrel{(\ref{real1})}{=} & \varphi^+_1(S^{-1}\R^{(1)} ) 
(a_2^*\tl S^{-1}\R^{(2)}) \nn
& \stackrel{(\ref{SR})}{=} & \varphi^+_1(\R^{(1)})
(a_2^*\tl \R^{(2)}) \nn
& \stackrel{(\ref{Def+})}{=} & \chi^+(a_2^*). \nonumber
\eea
Similarly one proves (\ref{result1}) for $\chi^-$.
Under the second assumptions, for any $a_2\in \c{A}_2$, 
\bea
[\chi^+(a_2) ]^* & \stackrel{(\ref{past1})}{=} & 
[(a_2\tl \R^{-1(2)})\, \varphi^+_1(S\R^{-1(1)} )]^* \nn
 &=& [\varphi^+_1(S\R^{-1(1)} ) ]^* (a_2\tl \R^{-1(2)})^* \nn
& \stackrel{(\ref{hypo2}),(\ref{starmod})}{=} & 
\varphi^-_1(S^{-1}\R^{-1(1)}{}^* ) 
(a_2^*\tl S^{-1}\R^{-1(2)}{}^*) \nn
& \stackrel{(\ref{real2})}{=} & \varphi^-_1(S^{-1}\R^{-1(2)} ) 
(a_2^*\tl S^{-1}\R^{-1(1)}) \nn
& \stackrel{(\ref{SR})}{=} & \varphi^-_1(\R^{-1(2)})
(a_2^*\tl \R^{-1(1)}) \nn
& \stackrel{(\ref{Def-})}{=} & \chi^-(a_2^*). \nonumber
\eea
By similar arguments one proves the claim for $\chi^-$. \ep

It should be noted that
there also exist non--standard star structures on $U_q \g$
for $|q|=1$, in particular the compact form ${X_i^{\pm}}^* = X_i^{\mp}$,
$K_i^* = K_i^{-1}$ in terms of the Cartan--Weyl generators. Then 
\be
\R^*=\R_{21}^{-1},                                       \label{real3}
\ee
while the coproduct 
does not fulfill (\ref{starhopf}) as in a standard Hopf
$*$-algebra but becomes flipped under the star.
This nevertheless has the correct classical limit,
because the coproduct is cocommutative for $q = 1$.
In certain cases (in particular on the fuzzy quantum sphere \cite{fuzzyq} 
discussed in section 6, but see also \cite{ads_space}), 
it is then possible to define a 
star structure on each $\c{A}_i$, which takes the form 
$a_{i;k}^* = \pm \Omega_i a_{i;k} \Omega_i^{-1}$ on the generators
$a_{i;k}$ of $\c{A}_i$. 
Here $\Omega_i = {}^4\!\!\sqrt{v_i}^{-1} \omega_i$, where
$v_i$ and $\omega_i$ are the realizations in $\c{A}_i$
(using an algebra map from $U_q \g$ to $\c{A}_i$ as above)
of the central element $v \in U_q\g$ (\ref{Defv}) and 
the ``universal Weyl element'' 
$\omega$ in an extension of $U_q\g$ 
\cite{kirill_resh}. 
All this must be defined 
in some representation of $\c{A}_i$; for more details 
see\footnote{The $v$ in \cite{fuzzyq,ads_space} is the square root of our 
$v$ here.} \cite{fuzzyq,ads_space}.
If moreover there exists an element $\Omega$  which realizes 
${}^4\!\!\sqrt{v}^{-1} \omega$ in $\c{A}^+=\c{A}_1\underline{\otimes}^+\c{A}_2$
or a ``physical'' subspace thereof, then 
it follows easily from (\ref{real3}) that the star structure 
$a_{i;k}^* = \pm \Omega a_{i;k} \Omega^{-1}$ on $\c{A}^+$
is consistent with the commutation 
relations of the braided tensor product algebra 
$\c{A}^+$. This star then has the
correct classical limit, and the same construction
also works for $ \c{A}^-$. 

\sect{Unbraiding `chains' of braided quantum Euclidean 
spaces or spheres}
\label{qspaces}

In this section we consider the braided tensor product of $M\ge 2$ copies
$\c{A}_1,\c{A}_2$, $...,\c{A}_M$ of the quantum Euclidean space $\b{R}_q^N$
~\cite{FadResTak89} (the $\uqs$-covariant quantum space),
i.e.  of the unital associative algebra generated by $x^i$ fulfilling the
relations 
\be 
\c{P}_a{}^{ij}_{hk}x^hx^k=0, \label{xxrel} 
\ee 
where $\c{P}_a$
denotes the $q$-deformed antisymmetric projector appearing in the decomposition
of the braid matrix $\hat R$ of $\uqs$ [given in formula (\ref{defRsoN})], 
or of the quotient space of the latter obtained by setting $r^2:=x^ix_i=1$
[the quantum $(N-1)$-dimensional sphere $S_q^{N\!-\!1}$].
The multiplet $(x^i)$ carries the fundamental vector representation $\rho$
of $\uqs$:  for any $g\in \uqs$
\be 
x^i\tl g=\rho^i_j(g)x^j.  \label{fund1} 
\ee
We shall enumerate the different copies of the quantum Euclidean space
or sphere by attaching an additional greek index to them, 
e.g.  $\alpha=1,2,...,M$.
The prescription (\ref{absbraiding'}) to glue 
$\c{A}_1,...,\c{A}_M$ into a $\uqs$-module 
associative algebra $\c{A}^-$ gives the following cross commutation relations
between their respective generators:
\be
x^{\alpha,i}x^{\beta,j}=\hat R^{ij}_{hk}x^{\alpha,h}x^{\beta,k}
                                                   \label{braid1}
\ee
whenever $\alpha<\beta$.
Note that prescriptions (\ref{absbraiding'}), (\ref{absbraiding}) go into
each other under the inverse reordering $1,2,...,M\to M,...,2,1$.
Applying iteratively Theorem \ref{maintheo} we shall be able to
completely unbraid this iterated tensor product.

To define $\varphi_1^{\pm}$ one actually needs a slightly enlarged version
\cite{CerFioMad00} of $\b{R}_q^N$ (or $S_q^{N\!-\!1}$). One has
to introduce some new generators $\sqrt{r_a}$, with
$1\le a\le \frac N 2$, together with their  
inverses $(\sqrt{r_a})^{-1}$, requiring that 
\be
r_a^2=\sum\limits_{h=-a}^a x^hx_h=\sum\limits_{h=-a}^a g_{hk}x^hx^k
\ee
(note that, having set
$n:=\left[\frac N 2\right]$, $r_n^2$ coincides with $r^2$).
Moreover for odd $N$ we add also 
also $\sqrt{x^0}$ and its inverse
as new generators). In fact, the commutation relations 
involving these new generators can be fixed consistently, and turn out to
be simply $q$-commutation relations.
$r$ plays the role of `deformed Euclidean distance' 
of the generic `point of coordinates' $(x^i)$ of $\b{R}_q^N$ from the `origin';
$r_a$ is the `projection' of $r$ on the `subspace' $x^i=0$, $|i|>a$. 
In the previous equation $g_{hk}$ denotes the `metric matrix' of 
$SO_q(N)$:
\be
g_{ij}=g^{ij}=q^{-\rho_i} \delta_{i,-j}.          \label{defgij}
\ee
It is a $SO_q(N)$-isotropic tensor and is a deformation of the 
ordinary Euclidean metric.
Here and in the sequel $n:=\left[\frac N 2\right]$ is the rank of 
$so(N)$, the indices take the values
$i=-n,\ldots,-1,0,1,\ldots n$ for $N$ odd,
and $i=-n,\ldots,-1, 1,\ldots n$ for $N$ even.
Moreover, we have introduced the notation
$(\rho_i)=(n-\frac{1}{2},\ldots,\frac{1}{2},0,-\frac{1}{2},
\ldots,\frac{1}{2}-n)$
for $N$ odd, $(n-1,\ldots,0,0,\ldots,1-n)$ for $N$ even. 
In the case of even $N$ one needs to include also
the FRT generator  $\c{L}^-{}^1_1$ and its inverse $\c{L}^+{}^1_1$
(which are generators of $U_qso(N)$ belonging to the Cartan 
subalgebra) among the generators of $\c{A}_1$. 
They satisfy the commutation relations
\be
\c{L}^-{}^1_1 x^{\pm 1}=q^{\pm 1} x^{\pm 1} \c{L}^-{}^1_1, \quad\quad
\c{L}^-{}^1_1 x^{\pm i}= x^{\pm i} \c{L}^-{}^1_1 \hbox{ for } i>1. 
\label{xkapparel}
\ee
with $\c{A}_1$, and the standard FRT relations with the rest
of $U_qso(N)$. One can easily show that the extension of the action of 
$U_qso(N)$
to $\sqrt{r_a},(\sqrt{r_a})^{-1}$ is uniquely determined by the constraints
the latter fulfil; it is a bit complicated and therefore will 
be omitted, since we will not need its explicit expression. 
The action of $H$ on $\c{L}^-{}^1_1$ is the standard (right) adjoint 
action. 
Note that the maps $\varphi^{\pm}_1$ have no analog
in the ``undeformed'' case ($q=1$),
because $\c{A}_1\equiv\b{R}^N$ is abelian, whereas 
$H\equiv U_qso(N)$ is not. 

The unbraiding procedure is recursive. We use the
homomorphism $\varphi_1$ found in Ref. \cite{CerFioMad00}
and start by unbraiding the first copy from the others. Following
Theorem \ref{maintheo}, we perform the following 
change of generators in $\c{A}^-$ 
\[
\ba{l}
y^{1,i}:=x^{1,i}                                       \\
y^{\alpha,i}:=\chi^-(x^{\alpha,i})
\stackrel{(\ref{Def-})}{=}\varphi_1(\R^{-1(2)})\rho^i_j
(\R^{-1(1)})x^{\alpha,j}=
\varphi_1(\c{L}^-{}^i_j)x^{\alpha,j}, \:\:\: \alpha>1.                     
\ea                                              
\]
In the last equality we have used the definition (\ref{frt})
of the FRT
generators \cite{FadResTak89} of $\uqs$.
In appendix \ref{appeuc} we recall the $\varphi^{\pm}$
images of the latter. In view of formula (\ref{imagel-})
we thus find
\be
\ba{l}
y^{1,i}:=x^{1,i}                                       \\
y^{\alpha,i}:=g^{ih}[\mu^1_h,x^{1,k}]_qg_{kj}\,
x^{\alpha,j}, \qquad \alpha>1.                     
\ea                                              \label{def1}
\ee
The suffix 1 in $\mu_a^1$ means that the special elements $\mu_a$
defined in (\ref{defmu}) must be taken as elements
of the first copy of $\b{R}^N_q$. In view of (\ref{defmu})
we see that $g^{ih}[\mu^1_h,x^{1,k}]_qg_{kj}$
are rather simple polynomials in $x^i$ and $r_a^{-1}$,
homogeneous of total degree 1 in the coordinates $x^i$ and $r_a$. 
Hence (\ref{def1}) is  a transformation
of polynomial type and therefore likely to be implemented as a
well-defined operator transformation also when representing
$\c{A}^-$ as an algebra of operators on some linear space.
Using the results (\ref{vf-}) given in the appendix
we give now the explicit expression of (\ref{def1}) for $N=3$:
\bea
&&y^{\alpha,-} = -qh\gamma_1 \frac r{x^0}x^{\alpha,-} \nn
&&y^{\alpha,0} = \sqrt{q}(q+1) \frac 1{x^0}x^+x^{\alpha,-}+x^{\alpha,0} \\
&&y^{\alpha,+} = \frac{\sqrt{q}(q+1)}{h\gamma_1 rx^0}(x^+)^2 x^{\alpha,-}+
\frac{q^{-1}+1}{h\gamma_1 r}x^+ x^{\alpha,0}- 
\frac 1{qh\gamma_1 r}x^0 x^{\alpha,+} \nonumber
\eea
for any $\alpha=2,...,M$. Here we have set $x^i\equiv x^{1,i}$, 
$h\equiv \sqrt{q}-1/\sqrt{q}$, replaced
for simplicity the values $-1,0,1$ of the indices by the ones $-,0,+$
and denoted  by $\gamma_1\in\b{C}$ a free parameter. 

As a consequence of the theorem we find
\begin{corollary}
\bea
&& [y^{1,i},y^{\alpha,j}]=0\qquad \alpha>1           \label{pippo1}\\
&& y^{\alpha,i}y^{\beta,j}=\hat R^{ij}_{hk}y^{\alpha,h}y^{\beta,k}
 \qquad 1<\alpha<\beta                                \label{pippo2}\\
&& \c{P}_a{}^{ij}_{hk}y^{\alpha,h}y^{\alpha,k}=0    \label{pippo3}
\eea
\end{corollary}

By (\ref{pippo1}) the subalgebra 
$\tilde\c{A}_1^-\equiv\c{A}_1$ of $\c{A}^-$ generated by
$y^{1,i}\equiv x^{1,i}$ commutes with the subalgebra generated
by $y^{2,i},...,y^{M,i}$, which we shall call $\tilde\c{A}^-$.
This was the first step of the unbraiding procedure. Now we can reiterate 
the latter for $\tilde\c{A}^-$, with $y^{2,i}$ playing the role
of $x^{1,i}$. After $M-1$ steps, we shall have determined $M$
independent commuting subalgebrae of $\c{A}^-$ which we
shall call $\tilde\c{A}_{\alpha}^-$, $\alpha=1,...,M$.

Theunbraiding procedure for the alternative braided tensor
product stemming
from prescription (\ref{absbraiding}) arises by iterating the
change of generators
\be
\ba{l}
y'{}^{M,i}:=x^{M,i}                                       \\
y'{}^{\alpha,i}:=\varphi_M(\c{L}^+{}^i_j)x^{\alpha,j}
=g^{ih}[\bar\mu^M_h,x^{M,k}]_{q^{-1}}g_{kj}\, x^{\alpha,j}, \qquad \alpha<M.
                                                   \label{def1'}
\ea
\ee
$\bar\mu_a^M$ are 
the special elements defined in (\ref{defbarmu}) belonging to
the $M$-th copy of $\b{R}_q^N$.
Using the results (\ref{vf+}) given in the appendix
we give the explicit expression of  (\ref{def1'}) for
$N=3$: for any $\alpha=1,...,M\!-\!1$,
\bea
&&y'{}^{\alpha,-} = -h\bar\gamma_1 \frac {z^0}{r_z}x^{\alpha,-}
+\frac{k\bar\gamma_1}{\sqrt{q}r_z}z^-x^{\alpha,0}
+\frac{q^{-2}k\bar\gamma_1}{r_zz^0}(z^-)^2x^{\alpha,+}\nn
&&y'{}^{\alpha,0} = x^{\alpha,0} +q^{-\frac 12}(q^{-1}+1) 
\frac 1{z^0}z^-x^{\alpha,+}\\
&&y'{}^{\alpha,+} = -
\frac {r_z}{h\bar\gamma_1 z^0}x^{\alpha,+} \nonumber
\eea
Here we have set $z^i\equiv x^{M,i}$,
$r_z^2\equiv x^{M,i}x^M_i$, $k\equiv q-q^{-1}$ and $\bar\gamma_1\in\b{C}$ 
is a free parameter.

Again, the subalgebra $\tilde\c{A}_M^+\approx \b{R}_q^N$ of $\c{A}^+$ 
generated by
$y^{M,i}\equiv x^{M,i}$ commutes with the subalgebra generated
by $y^{1,i},...,y^{M-1,i}$, which we shall call $\tilde\c{A}^+$.
This was the first step of the unbraiding procedure. Now we can reiterate 
the latter for $\tilde\c{A}^+$, with $y^{M\!-\!1,i}$ playing the role
of $x^{M,i}$. After $M-1$ steps, we shall have determined $M$
independent commuting subalgebrae of $\c{A}^+$ which we
shall call $\tilde\c{A}_{\alpha}^+$.

We summarize the results of this section:

\begin{prop} Let $\c{A}_1,\c{A}_2,...,\c{A}_M$ be $M$ 
copies of the \uqs-covariant quantum Euclidean space 
(or sphere).
Then $\c{A}_1\underline{\otimes}^{\pm}\c{A}_2\underline{\otimes}^{\pm}
...\underline{\otimes}^{\pm}\c{A}_M=\c{A}_1\otimes\tilde\c{A}_2^{\pm}
\otimes...\otimes\tilde\c{A}_M^{\pm}$, where 
$\tilde\c{A}_2^{\pm},...,\tilde\c{A}_M^{\pm}$ are subalgebras
of the lhs isomorphic to $\c{A}_1$.
\end{prop}
By a suitable choice of $\gamma_1,\bar\gamma_1$, as well as of the other
free parameters appearing in the definitions of $\varphi^{\pm}$
for $N>3$ (see appendix \ref{appeuc}), one can make $\varphi^{\pm}$
into $*$-homomorphisms
when $|q|=1$, and make them satisfy the relation
\be
[\varphi^{\pm}(g)]^*=\varphi^{\mp}(g^*)
\ee
when $q\in\b{R}^+$. Since these relations are of the type considered
in proposition \ref{propstar}, the claims of the latter for
$\chi^{\pm}$ and their consequences hold. In particular,
when $|q|=1$ one has a well-defined $*$ on the braided
tensor product of $\c{A}_1,...,\c{A}_M$ mapping each of the
independent, commuting subalgebras $\tilde\c{A}_i^{\pm}$ into itself.
On the contrary for real $q$ one can consider the map
$*:\c{A}^+\to\c{A}^-$ defined by (\ref{antimult}) or a
$*$-structure on $\c{A}^{\pm}$
defined in a way similar to what we have done in (\ref{altstar}),
\be
(x^{\alpha,i})^{\star}=x^{M\!-\alpha\!+\!1,j}g_{ji}.
\ee
The latter has not the
standard classical limit. A short calculation
shows that the latter implies
\be
(y^{\alpha,i})^{\star}=y'{}^{M\!-\alpha\!+\!1,j}g_{ji}.
\ee

\sect{Unbraiding `chains' of braided Heisenberg algebras}
\label{heisenberg}

In this section we consider the braided tensor product of $M\ge 2$ copies
$\c{A}_1,\c{A}_2$, $...,\c{A}_M$ of the 
$\uqg$-covariant deformed Heisenberg algebra 
$\c{D}_{\epsilon,\g}$, $\g=sl(N)$, $so(N)$ 
\cite{PusWor89,WesZum90,CarSchWat91}, 
i.e. the unital associative algebra generated by $x^i,\partial_j$
fulfilling the relations
\bea
&&\c{P}_a{}^{ij}_{hk}x^hx^k=0 \nn
&&\c{P}_a{}^{ij}_{hk}\partial_j\partial_i=0 \label{gringo}\\
&&\partial_i x^j=\delta^i_j+(q\gamma\hat R)^{\epsilon}
{}_{ih}^{jk}x^h\partial_k,
\nonumber
\eea
where  
$\gamma=q^{\frac 1N},1$ respectively for $\g=sl(N),so(\!N\!)$,
and the exponent $\epsilon$ can take either value $\epsilon=1,-1$.
$\hat R$ denotes the braid matrix of $\uqg$ [given in formulae
(\ref{defRslN}) and (\ref{defRsoN})], and again $\c{P}_a$ the
antisymmetric projector appearing in the decomposition
of the latter. 
The coordinates $x^i$ transform according to the fundamental
vector representation of $\uqg$, as in (\ref{fund1}), 
whereas the `partial derivatives' trasform according the
contragredient representation,
\be
\partial_i\tl g=
\partial_h\rho^h_i(S^{-1}g).                  \label{fund2}
\ee
The indices will take the values  $i=1,...,N$ if $\g=sl(N)$,
the same values considered in the previous section if $\g=so(N)$.
Clearly in the latter case $\c{D}_{\epsilon,\g}$ has the quantum
Euclidean space generated by $x^i$ as a module subalgebra. 

Again, we shall enumerate the different copies by attaching
to them an additional greek index, e.g. $\alpha=1,2,...,M$. 
The prescription (\ref{absbraiding}) to glue 
$\c{A}_1,...,\c{A}_M$ into a $\uqg$-module 
associative algebra $\c{A}^+$ (see also Ref. \cite{Fio98})
gives the following cross commutation relations 
between their respective generators 
\be
\ba{ll}
x^{\alpha,i}x^{\beta,j}=\hat R^{ij}_{hk}x^{\beta,h}x^{\alpha,k}\qquad \qquad&
\partial_{\alpha,i}\partial_{\beta,j}=
\hat R_{ji}^{kh}\partial_{\beta,h}\partial_{\alpha,k}\cr
\partial_{\alpha,i}x^{\beta,j}=
\hat R^{-1}{}^{jh}_{ik}x^{\beta,k}\partial_{\alpha,h}\qquad\qquad &
\partial_{\beta,i}x^{\alpha,j}=
\hat R^{jh}_{ik}x^{\alpha,k}\partial_{\beta,h}
\ea
\label{cross}
\ee
when $\alpha>\beta$. With respect to Ref. \cite{Fio98} we have
called the generators $x^i,\partial_j$ instead of $A^i,A^+_j$, inverted
the order of the product due to covariance w.r.t. the {\it right}
(instead of the {\it left}) \uqg-action,
and for the sake of simplicity we have put equal to one possible factors 
at the rhs of (\ref{cross}).

In Ref. \cite{Fiocmp95,ChuZum95}
algebra homomorphisms
$\varphi: \c{D}_{\epsilon,\g}\cocross H\rightarrow \c{D}_{\epsilon,\g}$ 
have been determined for $\g=so(n)$ and $\g=sl(N),so(N)$
respectively. This is the $q$-analog of vector field realization
of \g on the corresponding
$\g$-covariant (undeformed) space, e.g.
$\varphi_1(E^i_j)=x^i\partial_j-\frac 1N \delta^i_j$
in the $\g=sl(N)$ case. The searched maps
$\varphi^{\pm}$ will be simply the restrictions
of $\varphi$ to $\c{D}_{\epsilon,\g}\cocross H^{\pm}$. 
In Ref. \cite{Fiocmp95} there are among others
the $\varphi$-images of the Chevalley generators of 
\uqs\footnote{One should take care of the fact
that in Ref. \cite{Fiocmp95}
we considered \uqs acting by
a {\it left} action, instead of a right one, what
manifests itself in a replacement $q\to q^{-1}$, or equivalently
in an opposite coproduct. The rules for passing from right to left
are described in Sect. \ref{leftsec}.},
in Ref. \cite{ChuZum95} there are the $\varphi$-images of 
the generators of \uqg playing the role of ``vector fields'' on 
$G_q$. By the change of generators described in Ref. 
\cite{FadResTak89} 
one can easily pass from the
Chevalley to the FRT generators
$\c{L}^{\pm}{}^i_j$ (\ref{frt}), whereas the relation
between the latter and the vector fields is recalled in (\ref{defZ}). 
The FRT generators are the ones explicitly needed in writing down
$\chi^{\pm}(x^i)$ and $\chi^{\pm}(\partial_i)$.
For example, for $\g=sl(2)$ and $\epsilon=1$
one finds 
\bea
&&\varphi(\c{L}^+{}^1_1)=\varphi(\c{L}^-{}^2_2)=
[\varphi(\c{L}^-{}^1_1)]^{-1}=[\varphi(\c{L}^+{}^2_2)]^{-1}=
\alpha\Lambda^{\frac 12}\left[1\!+\!(q^2\!-\!1)
x^2\partial_2\right]^{\frac 12}\nn
&&\varphi(\c{L}^+{}^1_2)=-\alpha kq^{-1}\Lambda^{\frac 12}\left[1+(q^2-1)
x^2\partial_2\right]^{-\frac 12}x^1\partial_2 \\
&&\varphi(\c{L}^-{}^2_1)=\alpha kq^3\Lambda^{\frac 12}\left[1+(q^2-1)
x^2\partial_2\right]^{-\frac 12}x^2\partial_1,\nonumber
\eea
where $\alpha$ is fixed by (\ref{L+L-rel}) to be
$\alpha=\pm 1,\pm i$ and we have set
\be
\Lambda^{-2}:=1+(q^2-1)x^i\partial_i.        \label{defLambda-2'}
\ee
Whereas for $\g=so(3)$ and $\epsilon=1$
one finds on the positive Borel subalgebra
\be
\ba{l}
\varphi(\c{L}^+{}^-_-)=-\alpha\Lambda\left[1+(q-1)x^0\partial_0+(q^2-1)
x^+\partial_+\right]\\
\varphi(\c{L}^+{}^-_0)=\alpha k\Lambda(x^-\partial_0-\sqrt{q}x^0\partial_+)\\
\varphi(\c{L}^+{}^-_+)=\frac 1{1+q^{-1}}\varphi(\c{L}^+{}^-_0)
\varphi(\c{L}^+{}^0_+)\\
\varphi(\c{L}^+{}^0_0)=1 \\
\varphi(\c{L}^+{}^0_+)=-q^{-\frac 12}[\varphi(\c{L}^+{}^-_-)]^{-1}
\varphi(\c{L}^+{}^-_0)\\
\varphi(\c{L}^+{}^+_+)= [\varphi(\c{L}^+{}^-_-)]^{-1}
\ea \label{part1}
\ee
and on the negative Borel subalgebra
\be
\ba{l}
\varphi(\c{L}^-{}^-_-)=-\left(\alpha\Lambda\left[1+(q-1)x^0\partial_0+(q^2-1)
x^+\partial_+\right]\right)^{-1}\\
\varphi(\c{L}^-{}^0_-)=-\alpha q^2 k \varphi(\c{L}^-{}^-_-)
\Lambda(x^0\partial_--\sqrt{q}x^+\partial_0)\\
\varphi(\c{L}^-{}^+_-)=\frac1{1+q}\varphi(\c{L}^-{}^+_0)
\varphi(\c{L}^-{}^0_-)\\
\varphi(\c{L}^-{}^0_0)=1 \\
\varphi(\c{L}^-{}^+_0)=-\alpha q^{\frac 32}k\Lambda(x^0\partial_--
\sqrt{q}x^+\partial_0)\\
\varphi(\c{L}^-{}^+_+)= [\varphi(\c{L}^-{}^-_-)]^{-1}. 
\ea
\label{part2}
\ee
Here we have set
\be
\Lambda^{-2}:=[1+(q^2-1)x^i\partial_i+\frac{(q^2-1)^2}{\omega_1^2}
(g_{ij}x^ix^j)(g^{hk}\partial_k\partial_h)],        \label{defLambda-2}
\ee
where
$$
\omega_a:=(q^{\rho_a}+q^{-\rho_a}), 
$$
and replaced
for simplicity the values $-1,0,1$ of the indices by the ones $-,0,+$. 
In either case the $\varphi$-images of $\c{L}^+{}^i_j$ and
$\c{L}^-{}^j_i$ for $i>j$ vanish, because the latter do.

We see that strictly speaking $\varphi$ takes values in
some appropriate completion of $\c{D}_{\epsilon,\g}$, containing
at least the square root and inverse square root of the polynomial 
$\Lambda^{-2}$ respectively
 defined in (\ref{defLambda-2'}), (\ref{defLambda-2}),
as well as the square root of $[1+(q^2-1)
x^2\partial_2]$ and its inverse,
when $\g=sl(2)$, and the inverses
(\ref{part1})$_6$, (\ref{part2})$_6$, when $\g=so(3)$. 
Apart from this minimal completion, another possible one
is the socalled $h$-adic, namely the ring
of formal power series in $h=\log q$ with coefficients in
$\c{D}_{\epsilon,\g}$. Other completions, e.g. in operator norms,
can be considered according to the needs.
One can easily show that the extention of the action of $H$ 
to any such completion is uniquely determined (we omit to write 
down its explicit expression, since we don't  need it).

According to the main theorem, we set 
\be
\ba{rcl}
&&y^{1,i}\equiv x^{1,i}                                       \\
&&\partial_{y,1,a} \equiv  \partial_{1,a} \\
&&y^{\alpha,i} \equiv\chi^-(x^{\alpha,i})=
\varphi_1(\c{L}^-{}^i_j)x^{\alpha,j} \qquad \qquad \alpha>1\\
&&\partial_{y,\alpha,a}  \equiv\chi^-(\partial_{\alpha,a})=
\varphi_1(S\c{L}^-{}^d_a)\partial_{\alpha,d} \qquad \alpha>1 
\ea
\label{def2}
\ee
and we find
\begin{corollary}
\bea
&&\c{P}_a{}^{ij}_{hk}y^{\alpha,h}y^{\alpha,k}=0 \nn
&&\c{P}_a{}_{ij}^{hk}\partial_{y\alpha,k}
\partial_{y,\alpha,h}=0 \label{gigi}\\
&&\partial_{y,\alpha,i}y^{\alpha,j}=
\delta^j_i+(q\hat R)^{\epsilon_{\alpha}}{}^{jl}_{im}y^{\alpha,m}
\partial_{y,\alpha,l} \nonumber
\eea
for all $\alpha=1,...,M$, together with
\be
\ba{ll}
[y^{1,i},y^{\alpha,j}]=0\qquad\qquad &[\partial_{y,1,i}, y^{\alpha,j}]=0 \cr
[\partial_{y,\alpha,i}, y^{1,j}]=0 \qquad\qquad &
[\partial_{y,1,i},\partial_{y,\alpha,j}]=0
\ea
\label{giggi}
\ee
when $\alpha>1$, and
\be
\ba{l}
y^{\alpha,i}y^{\beta,j}=\hat R^{ij}_{hk}y^{\beta,h}y^{\alpha,k}\\
\partial_{y,\alpha,i}\partial_{y,\beta,j}=
\hat R_{ji}^{kh}\partial_{y,\beta,h}\partial_{y,\alpha,k}\\
\partial_{y,\alpha,i}y^{\beta,j}=
\hat R^{-1}{}^{jh}_{ik}y^{\beta,k}\partial_{y,\alpha,h}\\
\partial_{y,\beta,i}y^{\alpha,j}=
\hat R^{jh}_{ik}y^{\alpha,k}\partial_{y,\beta,h}
\ea
\label{ycross}
\ee
when $1<\beta<\alpha$.
\end{corollary}
By (\ref{giggi}) 
$y^{1,i}\equiv x^{1,i}$ and $\partial_{y,1,i}\equiv \partial_{1,i}$
commute with the subalgebra generated
by $y^{2,i},...,y^{M,i}$ and  $\partial_{y,2,i},...,\partial_{y,M,i}$
which we shall call $\tilde\c{A}^+$.
This was the first step of the unbraiding procedure. Now we can reiterate 
the latter for $\tilde\c{A}^+$, with $y^{2,i}, \partial_{y,2,i}$ 
playing the role
of $x^{1,i},\partial_{1,i}$. After $M-1$ steps, we shall have determined $M$
independent commuting subalgebras of $\c{A}^+$ which we
shall call $\tilde\c{A}_{\alpha}^+$.

For the sake of brevity we omit
the unbraiding procedure for the alternative braided tensor
product algebra stemming
from prescription (\ref{absbraiding'}), which can be found following
arguments completely analogous to the ones presented at the
end of Section \ref{qspaces}. We summarize the results of this section by

\begin{prop} Let $\c{A}_1,\c{A}_2,...,\c{A}_M$ be $M$ 
copies of the \uqg-covariant
Hei\-senberg algebra $\c{D}_{\epsilon,\g}$, $\g=sl(N)$, $so(N)$.
Then $\c{A}_1\underline{\otimes}^{\pm}\c{A}_2\underline{\otimes}^{\pm}
...\underline{\otimes}^{\pm}\c{A}_M=\c{A}_1\otimes\tilde\c{A}_2^{\pm}
\otimes...\otimes\tilde\c{A}_M^{\pm}$, where 
$\tilde\c{A}_2^{\pm},...,\tilde\c{A}_M^{\pm}$ are subalgebras
of the lhs isomorphic to $\c{D}_{\epsilon,\g}$.
\end{prop}

Relations (\ref{|q|=1}), (\ref{starmod}), (\ref{fund1}), 
(\ref{fund2}) and (\ref{gringo}) fix
the $*$-structure of $\c{A}_1$ to be
\be
(x^i)^*=x^i, \qquad (\partial_i)^*=-\partial_i\cases{q^{\pm 2(N-i+1)}
\:\:\mbox{ if }H=U_qsl(N)\cr
q^{\pm N+\rho_i}\:\:\mbox{ if }H=U_qso(N)}
\ee
if $|q|=1$, and 
\be
(x^h)^*=x^kg_{kh}, \qquad (\partial_i)^*=-\frac{\Lambda^{\pm 2}}{
q^{\pm N}+q^{\pm 2}}\left[(g^{jh}\partial_h\partial_j),\, x^i\right]
\ee
if $H=U_qso(N)$ and $q\in\b{R}^+$. The upper or lower sign respectively refer
to the choices
$\epsilon=1,-1$ in (\ref{gringo})$_3$, and 
$\Lambda^{\pm 2}$ are respectively defined by
\be
\Lambda^{\pm 2}:=\left[1+(q^{\pm 2}-1)x^i\partial_i+ \frac{(q^{\pm 2}
-1)^2}{\omega_n^2} r^2 (g^{ji}\partial_i\partial_j)\right]^{-1}.
\ee
The map $\varphi$ is a $*$-homomorphism both for
$q$ real and $|q|=1$. If we denote by
$\varphi^{\pm}$ its restrictions to 
$\c{A}\cocross H^{\pm}$, then they are $*$-homomorphisms
when $|q|=1$ (see appendix \ref{appstar}), and 
fulfill the relation
\be
[\varphi^{\pm}(g)]^*=\varphi^{\mp}(g^*)
\ee
when $q\in\b{R}^+$ \cite{Fiocmp95}. 
Since these relations are of the type considered
in proposition \ref{propstar}, the claims of the latter for
$\chi^{\pm}$ and their consequences hold. In particular,
when $|q|=1$ one has a well-defined $*$ on the braided
tensor product of $\c{A}_1,...,\c{A}_M$ mapping each of the
independent, commuting subalgebras 
$\tilde\c{A}_{\alpha}^{\pm}$ into itself.

\medskip

Finally, the above results have an important corollary.
According to Hochschild cohomology arguments developed by
Gerstenhaber \cite{Ger64} and applicable to Heisenberg algebras
because of the results found by Du Cloux 
in Ref. \cite{Duc85}, any deformed
Heisenberg algebra, in particular the braided tensor products
of $\c{A}_1,...,\c{A}_M$ considered in this section, can be
realized simply by a change of generators in the $h$-adic
completion, $h=\log q$, of its undeformed counterpart (but in general 
{\it not} in other, e.g. {\it operator-norm}, completions). 
However explicit realizations are not provided by
these results. The results presented here, combined to some
older ones, allow to determine one such realization.
In Ref. \cite{Ogi92} Ogievetsky found an explicit realization $\phi$ or
`deforming map' of the elements of $\c{D}_{\epsilon,\g}$
in terms of formal power series in $h=\log q$ with coefficients in
the corresponding undeformed Heisenberg algebra. Another, less
explicit, one was found in Ref. \cite{Fiormp00}.
The composition of the unbraiding map 
found in this section, which allows to `decouple'
$M$ different copies of $\c{D}_{\epsilon,\g}$ from each other, 
with the map $\phi$ provides an explicit realization 
or `deforming map' of
the larger Heisenberg algebra $\c{A}$ (what we have
called the `braided chain of Heisenberg algebras'), 
in the $h$-adic completion of the
undeformed $(N\cdot M)$-dimensional Heisenberg algebra. 

\sect{Formulae for the left action}
\label{leftsec}

For psychological reasons we often prefer to work with
a left action rather than with a right one.
In this section we give the analogs for left $H$-module algebras
of the main results found so far for right $H$-module algebras.
The left action of $g\in H$ on a product fulfills
\bea
&&(gg')\tr a =g\tr (g'\tr a),                     \label{modalg1l}\\
&&g\tr (aa') =(g_{(1)}\tr a)(g_{(2)}\tr a').      \label{modalg2l}
\eea
The product laws in the braided tensor product  algebras
$\hat\c{A}^+,\hat\c{A}^-$ are respectively given by
\bea
&& a_2a_1= (\R^{-1(1)}\tr a_1)\, (\R^{-1(2)}\tr a_2). \label{absbraidingl'}\\
&& a_2a_1= (\R^{(2)}\tr a_1)\, (\R^{(1)}\tr a_2),     \label{absbraidingl}
\eea
The analog of Theorem \ref{maintheo} reads
\begin{theorem} 
Let $\{H,\R\}$ be a quasitriangular Hopf algebra and 
$H^+,H^-$ be Hopf subalgebras of $H$ such that $\R\in H^+\otimes H^-$.
Let $\hat\c{A}_1, \hat\c{A}_2$ be respectively a (left) $H^+$- and
a $H^-$-module algebra, so that we can define $\hat\c{A}^+$ 
as in (\ref{absbraidingl'}), and 
$\hat\varphi_1^+:H^+\cross \hat\c{A}_1\to \hat\c{A}_1$ be an algebra 
homomorphism fulfilling (\ref{ident0}), so that we can define
a map $\hat\chi^+:\hat\c{A}_2\rightarrow \hat\c{A}^+$ by
\be
\hat\chi^+(a_2):= (\R^{(2)}\tr a_2)\,\hat\varphi_1^+(\R^{(1)}).
\label{Def+l}
\ee
Alternatively, let $\hat\c{A}_1,\hat\c{A}_2$ be respectively a (left)
$H^-$- and a $H^+$-module algebra, so that we can define $\hat\c{A}^-$ 
as in (\ref{absbraidingl}), and 
$\hat\varphi_1^+:H^+\cross \hat\c{A}_1\to \hat\c{A}_1$ be an algebra 
homomorphism fulfilling (\ref{ident0}), so that we can define
a map $\hat\chi^-:\hat\c{A}_2\rightarrow \hat\c{A}^-$ by
\be
\hat\chi^-(a_2):= (\R^{-1(1)}\tr a_2)\,\hat\varphi_1^-(\R^{-1(2)}).
\label{Def-l}
\ee
In either case $\hat\chi^{\pm}$ are then injective algebra 
homomorphisms and 
\be
[\hat\chi^{\pm}(a_2),\hat\c{A}_1]=0,                       \label{commul}
\ee
namely the subalgebras 
$\tilde{\hat\c{A}}_2^{\pm}:=\hat\chi^{\pm}(\hat\c{A}_2)
\approx\hat\c{A}_2$ commute with $\hat\c{A}_1$. Moreover 
$\hat\c{A}^{\pm}=\hat\c{A}_1\otimes\tilde{\hat\c{A}}_2^{\pm}$.
\end{theorem}

The results of section \ref{star} apply without modifications
(one just has to place a $\hat{}$ in the appropriate places).

To enumerate the generators of the
algebras considered in Sections \ref{qspaces},\ref{heisenberg}
we shall exchange lower with upper indices,  so the generators
will read $x_{\alpha,i},\partial^{\alpha,i}$. This is necessary
if we wish the $x$'s to carry what we shall consider the
fundamental (vector) representation $\rho$ of \uqg, 
\be
g\tr x_i=x_j\rho_i^j(g),
\ee
rather than its contragredient $\rho^T\circ S$, because this
follows from the row$\times$col\-umn multiplication law
$\rho^i_h(gg')=\rho^i_j(g)\rho^j_h(g')$. Apart from this
replacement, all the commutation relations remain the same,
but can be rephrased in an equivalent way
exchanging lower with upper indices also
in the braid matrices and in the projectors $\c{P}_a$, because
$\hat R^T=\hat R$, $\c{P}_a{}^T=\c{P}_a$. For instance, the 
analog of (\ref{xxrel})
will read
\be 
\c{P}_a{}_{ij}^{hk}x_hx_k=0.                \label{xxrell} 
\ee 
The analogs of (\ref{fund1}), (\ref{fund2}) read
\bea
&& g\tr x_i=\rho_i^j(g)x_j  \label{fund1l} \\
&& g\tr \partial^i =
\partial^h\rho^i_h(Sg).                  \label{fund2l}
\eea

Algebra homomorphisms $\hat\varphi_1^{\pm}$ for the
algebras considered in Sections \ref{qspaces},\ref{heisenberg} are
immediately obtained in terms of the $\varphi_1^{\pm}$
described there, according to the rule
\be
\hat\varphi_1^{\pm}(\c{L}^{\pm}{}^h_j)=
U^{-1}{}^j_a\varphi_1^{\mp}(\c{L}^{\mp}{}^a_b) U^b_h.
\ee
Here 
\be
U^b_c:=\rho^b_c(u),
\ee 
$u\in H$ is a special
element as in (\ref{interantip}
), and at the rhs
the correct expression in the new notation
has lower and upper indices exchanged. 
If $\hat\c{A}_1$ is the quantum Euclidean space $\b{R}_q^N$
one finds, for instance,
\be
\hat\varphi_1^-(\c{L}^-{}^h_j)=
U^{-1}{}^j_a\,g_{ac}[\bar\mu^c,x_k]_{q^{-1}} g^{kb} U^h_b
\stackrel{(\ref{defU})}{=}g_{cj}[\bar\mu^c,x_k]_q g^{hk},
\ee
where $\mu^c$ is the same as $\mu_c$ [see \ref{defmu})],
but in the new notation.
For instance, when $|c|>1$ it reads 
\be
\bar\mu^c=\bar\gamma_c r_{|c|}^{-1}r_{|c|-1}^{-1} x_{-c},
\ee
with $\gamma_c$ defined as in (\ref{bargamma}) and $r_a$
($a\ge 0$) defined by the condition
$$
r_a^2=\sum\limits_{h=-a}^a x_hx^h=\sum\limits_{h=-a}^a g^{hk}x_hx_k.
$$
The analog of (\ref{def1}) is therefore (with $\alpha>1$)
\bea
&&y_{1,i}:=x_{1,i}                                   \label{def1l}    \\
&&y^{\alpha,i}:=\hat\chi^-(x_{\alpha,i})=x_{\alpha,j}\hat
\varphi_1(\c{L}^-{}^j_i)=x_{\alpha,j}
g_{hi}[\bar\mu^{1,h},x_{1,k}]_{q^{-1}}g^{jk}.        
\eea

\sect{Unbraiding `chains' of fuzzy quantum spheres}
\label{q_fuzzy}

As a last example, we consider the braided tensor product of $M$ copies
$\c{A}_1,...,\c{A}_M$ of the $q$--deformed fuzzy sphere $\hat S^2_{q,N}$
~\cite{fuzzyq}\footnote{To relate this to our conventions, 
the $q$ in \cite{fuzzyq}
should be replaced by $q^{-1/2}$}, which we consider 
as a left $U_q so(3)$ module algebra. It is generated by $x_i$ fulfilling the
relations 
\bea
\varepsilon_k^{i j}  x_i  x_j &=& \Lambda_N \; x_k, \nonumber\\
g^{i j} x_i x_j &=&  R^2.
\label{xxrel_q} 
\eea
Here $R>0$, 
\be
C_N = \frac{[N]_q [N+2]_q}{[2]_q^2}, \qquad
\Lambda_N = R \; \frac{[2]_{q^{N+1}}}{\sqrt{[N]_q [N+2]_q}}
\ee
where
$[n]_q:= \frac{q^{n/2}-q^{-n/2}}{q^{1/2}-q^{-1/2}}$, and
\be
\begin{array}{ll} 
\varepsilon^{1 0}_{1} = {q}^{1/2}, & \varepsilon^{0 1}_{1} = -q^{-1/2},\\
\varepsilon^{0 0}_{0} = q^{1/2}-q^{-1/2}, 
  & \varepsilon^{1 -1}_{0} = 1 = -\varepsilon^{-1 1}_{0}, \\
\varepsilon^{0 -1}_{-1} = q^{1/2}, & \varepsilon^{-1 0}_{-1} = -q^{-1/2}
\end{array}
\label{C_ijk}
\ee
are the spin 1 Clebsch--Gordan coefficients.
The multiplet $(x_i)$ carries the fundamental vector representation $\rho$
of $H = U_q so(3)$:  
\be 
g \tr x_i =x_j\rho^j_i(g).  
\ee
There is no obvious generalization to 
higher dimensions, but this algebra appears to be relevant e.g. to 
$D$--branes on the $SU(2)$ WZW model \cite{alekseev}. 
It has a unique irreducible 
representation, which is equivalent to $Mat(N+1)$.
Here we only consider the case 
$q \in \b{R}^+$, where the star structure is given by 
$x_i^* = g^{ij} x_j$. Then $\hat S^2_{q,N}$ is simply the 
``discrete series'' of 
Podles's spheres \cite{podles}. It was shown in  \cite{fuzzyq} that there
is a star--algebra homomorphism 
$\hat\varphi: H \cross \hat S^2_{q,N}\to \hat S^2_{q,N}$, which
takes a particularly simple form
\bea
\hat\varphi(E^+) &=& \frac 1{R} \sqrt{q^{-1}[2]_q C_N} \; x_1, \quad
\hat\varphi(E^-) = -\frac{1}{R} \sqrt{q [2]_q C_N}\;  x_{-1}, \nonumber\\
\hat\varphi(q^{H/2}) &=& \frac{[2]_{q^{N+1}}}{[2]_q}
         \Big(1- \frac{q^{1/2}-q^{-1/2}}{\Lambda_N}\;x_0 \Big)
\label{U_A_relation}
\eea
where $E^{\pm} = X^{\pm} q^{H/4} \in U_q so(3)$. Note that
$(1- \frac{q^{1/2}-q^{-1/2}}{\Lambda_N}\;x_0)$ is invertible 
since the eigenvalues of $q^{H/2}$ are positive (assuming $q>0$),
therefore $\hat\varphi(q^{-H/2}) \in \hat S^2_{q,N}$ is well--defined 
also. Hence the algebra homomorphisms $\hat \varphi$ is 
defined on the entire algebra $U_q so(3)$.
Using the definition (\ref{frt}) and the explicit form
for the universal $\R$ (see e.g. \cite{CH_P}), one finds 
\be
[\c{L}^-{}^i_j]=\left[\ba{ccc}
q^{-H/2}, \:&0, &0 \\
-(1-q^{-1})\sqrt{[2]_q} E^-, \:& 1, & 0 \\
q^{-1/2}(1\!-\!q^{-1})^2 q^{-H/2} (E^-)^2,\;& -(1-q^{-1})\sqrt{[2]_q}\; 
 q^{-H/2}E^-, & q^{H/2} 
\ea\right]
\label{vffuzzy-} 
\ee
and
\be
[\c{L}^+{}^i_j]\!=\!\left[\ba{ccc}
q^{-H/2},\:\: & (q-1)\sqrt{[2]_q}\; q^{-H/2}E^+,\: & 
(q-1)^2 q^{-H/2} (E^+)^2,\;\\
0, & 1, & q^{-1/2}(q-1)\sqrt{[2]_q} E^+ \\
0,& 0,& q^{H/2}
\ea\right].
\label{vffuzzy+} 
\ee
The unbraiding procedure then works as in Theorem 2. To be specific, assume
that the braided tensor product algebra is as in (\ref{absbraidingl'}). Then 
we set 
\bea
&&y_{1,i}:=x_{1,i}                                       \\
&&y_{\alpha,i}:=\hat\chi(x_{\alpha,i})=
  x_{\alpha,j}\hat\varphi_1(\c{L}^+{}^j_i), \qquad \alpha>1,       
\eea
without spelling out these expressions further. 
According to Theorem 2, they satisfy
\begin{corollary}
\bea
&&\varepsilon_k^{i j}  y_{\alpha,i}\;  y_{\alpha,j} = 
           \Lambda_N \; y_{\alpha,k}, \nonumber\\
&&g^{i j}  y_{\alpha,i} \; y_{\alpha,j} =  R^2   \nonumber 
\eea
for all $\alpha=1,...,M$, together with
\bea
&&[y_{1,i}, y_{\alpha,j}]=0 \\
&&y_{\alpha,i} y_{\beta,j}=\hat R_{ij}^{hk}y_{\beta,h} y_{\alpha,k} 
\label{ycross_fuzzy}
\eea
when $1<\alpha$ and $\alpha\beta$.
\end{corollary}
Iterating this procedure as before, we find 
 
\begin{prop} Let $\c{A}_1,\c{A}_2,...,\c{A}_M$ be $M$ 
copies of the $U_qso(3)$--covariant fuzzy quantum sphere.
Then $\c{A}_1\underline{\otimes}^{\pm}\c{A}_2\underline{\otimes}^{\pm}
...\underline{\otimes}^{\pm}\c{A}_M=\c{A}_1\otimes\tilde\c{A}_2^{\pm}
\otimes...\otimes\tilde\c{A}_M^{\pm}$, where 
$\tilde\c{A}_2^{\pm},...,\tilde\c{A}_M^{\pm}$ are subalgebras
of the lhs isomorphic to $\c{A}_1$.
\end{prop}

\app{Appendix}
\renewcommand{\theequation}{\thesubsection.\arabic{equation}}

\subsection{The universal $R$-matrix}                       \label{UnivR} 

In this appendix we recall the basics about the universal $R$-matrix
\cite{Dri86} of the quantum groups $\uqg$, while fixing our conventions.
Recall that the
universal $R$-matrix $\R$ is a special element
\be
\R\equiv\R^{(1)}\otimes\R^{(2)}\in\uqg\otimes \uqg
\label{swee}
\ee
intertwining between  $\Delta$ and opposite coproduct
$\Delta^{op}$, and so does also $\R^{-1}_{21}$:
\be
\ba{l}
\R( g_{(1)}\otimes g_{(2)})=(g_{(2)}\otimes g_{(1)})\R,  \\
\R^{-1}_{21}(g_{(1)}\otimes g_{(2)})=(g_{(2)}\otimes g_{(1)})\R^{-1}_{21}.  
\ea
\label{inter}
\ee
In  (\ref{swee}) we have used a Sweedler notation
with upper
indices: the right-hand side is a short-hand notation
for a sum $\sum_I\R_I^{(1)}\otimes\R_I^{(2)}$ of infinitely 
many terms.
We recall some useful formulae 
\bea
&&(\Delta \otimes \mbox{id})\R=\R_{13}\R_{23}, \label{delta1}\\
&&(\mbox{id}\otimes\Delta  )\R=\R_{13}\R_{12}, \label{delta2}\\
&&(S \otimes \mbox{id})\R=\R^{-1}= (\mbox{id}\otimes S^{-1})\R,
                                               \label{SR} \\
&& S^{-1}(g) = u^{-1}S(g) u.                  \label{interantip}
\eea
Here $u$ is any of the elements $u_1,u_2,..u_8$ defined below:
\be
\begin{array}{ll}
u_1:= (S\R^{(2)}) \R^{(1)} \qquad\qquad &u_2:= (S\R^{-1}{}^{(1)})
\R^{-1}{}^{(2)} \cr
u_3:= \R^{(2)} S^{-1}\R^{(1)} \qquad\qquad &u_4:= \R^{-1}{}^{(1)} 
S^{-1}\R^{-1}{}^{(2)} \cr
(u_5)^{-1}:= \R^{(1)} S\R^{(2)} \qquad\qquad &(u_6)^{-1}:= 
(S^{-1}\R^{(1)}) \R^{(2)} \cr
(u_7)^{-1}:= \R^{-1}{}^{(2)} S\R^{-1}{}^{(1)} 
\qquad\qquad &(u_8)^{-1}:= 
(S^{-1}\R^{-1}{}^{(2)}) \R^{-1}{}^{(1)} 
\label{defu}
\end{array}
\ee
In fact, using the results of Drinfel'd \cite{Dri86,Dri90} one can show that
\be
u_1=u_3=u_7=u_8=vu_2=vu_4=vu_5=vu_6,                      \label{Defv}
\ee
where $v$ is a suitable element belonging to the center of \uqs.

\medskip

{From} (\ref{inter}) and (\ref{delta1},\ref{delta2}) it follows the
universal Yang-Baxter relation
\be
\R_{12}\R_{13}\R_{23}=\R_{23}\R_{13}\R_{12},          \label{YBEQ1}
\ee
whence the other two relations follow
\bea
\R^{-1}{}_{12}\R^{-1}{}_{13}\R^{-1}{}_{23}& =& 
\R^{-1}{}_{23}\R^{-1}{}_{13}\R^{-1}{}_{12},\label{YBEQ2}\\
\R_{13}\R_{23}\R^{-1}{}_{12}& =& 
\R^{-1}{}_{12}\R_{23}\R_{13}.\label{YBEQ3}
\eea
As before, let $\rho$ be the fundamental $N$-dimensional
representation of $\g=sl(N),so(N),sp(N)$
By applying $\mbox{id}\otimes \rho^a_c\otimes\rho^b_d$ to 
(\ref{YBEQ1}), $\rho^a_c\otimes\rho^b_d\otimes\mbox{id}$ to
(\ref{YBEQ2}) and $\rho^a_c\otimes\mbox{id}\otimes\rho^b_d$
to (\ref{YBEQ3}) we respectively find the commutation relations
\bea
&&\hat R^{ab}_{cd}\,\c{L}^+{}^d_f\c{L}^+{}^c_e=
\c{L}^+{}^b_c\c{L}^+{}^a_d\,\hat R^{dc}_{ef} ,     \label{L+L+rel}\\
&&\hat R^{ab}_{cd}\,\c{L}^-{}^d_f\c{L}^-{}^c_e=
\c{L}^-{}^b_c\c{L}^-{}^a_d\,\hat R^{dc}_{ef}, \label{L-L-rel}\\
&&\hat R^{ab}_{cd}\,\c{L}^+{}^d_f\c{L}^-{}^c_e=
\c{L}^-{}^b_c\c{L}^+{}^a_d\,\hat R^{dc}_{ef},\label{L+L-rel}
\eea
where $\c{L}^{\pm}{}_l^a$ are the Faddeev-Reshetikin-Takhtadjan
generators \cite{FadResTak89} of $\uqg$, defined by
\be
\c{L}^+{}_l^a:=\R^{(1)}\rho_l^a(\R^{(2)})\qquad\qquad
\c{L}^-{}_l^a:=\rho_l^a(\R^{-1}{}^{(1)})\R^{-1}{}^{(2)}. \label{frt}
\ee

It is known \cite{FadResTak89} that $\{\c{L}^+{}^i_j,\c{L}^-{}^i_j\}$ and the
square roots of the
elements $\c{L}^{\pm}{}^i_i$  provide a (overcomplete)
set of generators of \uqg. 
Since in our conventions  
\be
\R\in H^+\otimes H^-,
\ee
then
$\c{L}^+{}_l^a\in H^+$ and $\c{L}^-{}_l^a\in H^-$. 
Beside (\ref{L+L+rel}-\ref{L+L-rel}) 
these generators fulfill 
\bea
&&\c{L}^+{}^i_j=0,           \hspace{1.5cm}\mbox{if $i>j$}\label{sfilza1}\\
&&\c{L}^-{}^i_j=0,           \hspace{1.5cm}\mbox{if $i<j$}\label{sfilza2}\\
&&\c{L}^{-}{}^i_i\c{L}^{+}{}^i_i=\c{L}^{+}{}^i_i\c{L}^{-}{}^i_i=1,
\hspace{1cm}\forall i  \label{sfilza3}\\
&&\c{L}^{\pm}{}^{-n}_{-n}....\c{L}^{\pm}{}^n_n=1,         \label{sfilza4}
\eea
and, when $\g=so(N),sp(N)$, some additional relations.
When $\g=so(N)$ the latter read
\be
\c{L}^{\pm}{}^i_j\c{L}^{\pm}{}^h_k g^{kj}=g^{hi}\qquad 
\c{L}^{\pm}{}_i^j\c{L}^{\pm}{}_h^k g_{kj}=g_{hi},  \label{LLg}
\ee
where $g_{ij}$ has been defined in (\ref{defgij}).
The braid matrix $\hat R$ is related to $\R$ by
$\hat R^{ij}_{hk}\equiv R^{ji}_{hk}:=(\rho^j_h\otimes\rho^i_k)\R$.
With the indices' convention described in sections \ref{qspaces},
\ref{heisenberg} $\hat R$ is given by
\be
\hat R = q^{-\frac 1N}\left[q \sum_i e^i_i \otimes e^i_i +
\sum_{\scriptstyle i \neq j} e^j_i \otimes e^i_j
+k \sum_{i<j} e^i_i \otimes e^j_j \right]            \label{defRslN}
\ee
when $\g=sl(N)$, and by
\bea
\hat R&=&q \sum_{i \neq 0} e^i_i \otimes e^i_i +
\sum_{\stackrel{\scriptstyle i \neq j,-j} 
{\mbox{ or } i=j=0}} e^j_i \otimes e^i_j+ q^{-1} 
\sum_{i \neq 0} e^{-i}_i
\otimes e^i_{-i} 
\label{defRsoN} \\
&&+k (\sum_{i<j} e^i_i \otimes e^j_j- 
\sum_{i<j} q^{-\rho_i+\rho_j} 
e^{-j}_i \otimes e^j_{-i}) \nonumber
\eea
when $\g=so(N)$.
Here $e^i_j$ is the $N \times N$ matrix with all elements
equal to zero except for a $1$ in the $i$th column and $j$th row.
The braid matrix of $sl(N)$ admits the orthogonal projector
decomposition
\be
\hat R = q\c{P}_S - q^{-1}\c{P}_a, \qquad\qquad \g=sl(N);
\label{projectorR}  
\ee
$\c{P}_a,\c{P}_S$ are the $U_qsl(N)$-covariant
deformed antisymmetric and symmetric projectors. 
The braid matrix of $so(N)$ admits the orthogonal projector
decomposition
\be
\hat R = q\c{P}_s - q^{-1}\c{P}_a + q^{1-N}\c{P}_t\qquad\qquad \g=so(N);      
\label{projectorR'}
\ee
$\c{P}_a,\c{P}_t,\c{P}_s$ are the $q$-deformed antisymmetric, trace,
trace-free symmetric projectors. 

The compact section of \uqg requires $q\in\b{R}^+$ if $\g=so(N)$, 
$q\in\b{R}$ if $\g=sl(N)$ and is characterized by the $*$-structure
\be
(\c{L}^{\pm}{}^i_j)^*=S\c{L}^{\mp}{}^j_i.                 \label{qReal}
\ee
For $\g=so(N)$ this amounts to
\be
(\c{L}^{\pm}{}^i_j)^*=g_{ih}\c{L}^{\mp}{}^h_k g^{kj}.        \label{qreal0}
\ee
The non-compact sections of \uqg require $|q|=1$ and are
characterized by the $*$-structure
\be
(\c{L}^{\pm}{}^i_j)^*=U^{-1}{}^i_r\,\c{L}^{\pm}{}^r_s\, U^s_j= 
u\,\c{L}^{\pm}{}^i_ju^{-1}.                                \label{|q|=1}
\ee
This can be checked using the property
$(\hat R^{ij}_{hk})^*=\hat R^{-1}{}^{ji}_{kh}$. Here we have defined
\be
U^i_j=\rho^i_j(u)
\ee
with $u$ any of the elements defined in (\ref{defu}). For $\g=so(N)$
one can take
\be
U^i_j:=g^{ih}g_{jh}.                                \label{defU}
\ee

{From} formulae (\ref{delta1}), (\ref{delta2}) in the 
Appendix~\ref{UnivR} one finds that the coproducts are given by
\be
\Delta(\c{L}^+{}^i_j)=\c{L}^+{}^i_h\otimes\c{L}^+{}^h_j \qquad\qquad
\Delta(\c{L}^-{}^i_j)=\c{L}^-{}^i_h\otimes\c{L}^-{}^h_j.
\label{coprodL}
\ee

\subsect{Proof of Proposition \ref{propalt}} 
We make use of the identity
\be
\varphi^{\pm} (g^{\pm}) \tl h^{\pm} =\varphi^{\pm} (g^{\pm} \tl h^{\pm}), 
                                                 \label{int}
\ee
for any $g^{\pm},h^{\pm} \in H^{\pm} $, which we prove 
in Ref. \cite{Fio00}. 
The right action appearing at the rhs is the (right) adjoint action
on itself
\be
h\tl g=Sg_{(1)} h g_{(2)},          \qquad\qquad g,h\in H; \label{adjo}
\ee
where $S$ denotes the antipode of the Hopf algebra $H$. 
We shall also need the inverse of (\ref{absbraiding}),
\be
a_1a_2= (a_2\tl \R^{-1(2)})\, (a_1\tl \R^{-1(1)}).     \label{invbraiding}
\ee
Now,
\bea
\chi^+(a_2) & \stackrel{(\ref{Def+})}{=}& 
\varphi_1^+(\R^{(1)})\, (a_2\tl \R^{(2)}) \nn
& \stackrel{(\ref{invbraiding})}{=}&  (a_2\tl \R^{(2)}\R^{-1(2')})
\, [\varphi_1^+(\R^{(1)})\tl \R^{-1(1')}]\nn
& \stackrel{(\ref{int})}{=}& (a_2\tl \R^{(2)}\R^{-1(2')})
\, \varphi_1^+(\R^{(1)}\tl \R^{-1(1')})\nn
& \stackrel{(\ref{adjo})}{=}& (a_2\tl \R^{(2)}\R^{-1(2')})
\, \varphi_1^+(S\R^{-1(1')}_{(1)}\,\R^{(1)}\, \R^{-1(1')}_{(2)})\nn
& \stackrel{(\ref{delta1})}{=}& (a_2\tl \R^{(2)}\R^{-1(2')}\R^{-1(2")})
\, \varphi_1^+(S\R^{-1(1")}\,\R^{(1)}\R^{-1(1')})\nn
& =& (a_2\tl \R^{-1(2")})
\, \varphi_1^+(S\R^{-1(1")}), \nonumber
\eea
which proves (\ref{past1}). Similarly one proves (\ref{past2}).

\subsect{The maps $\varphi^{\pm}$ for the quantum Euclidean spa\-ces
or spheres}
\label{appeuc}

We introduce the short-hand notation
\be
[A,B]_x=AB-xBA.                             \label{short-hand}
\ee
In Ref. \cite{CerFioMad00} we have found algebra homomorphisms
$\varphi^{\pm}:\b{R}_q^N\cocross U_q^{\pm}so(N)\to \b{R}_q^N$.
The images of
$\varphi^-$ on the negative FRT generators
read
\be
\varphi^-(\c{L}^-{}^i_j)=g^{ih}[\mu_h,x^k]_qg_{kj},      \label{imagel-}
\ee
where
\be
\begin{array}{ll}
\mu_0=\gamma_0 (x^0)^{-1}
&\quad\mbox{for $N$ odd,} \\[6pt]
\mu_{\pm 1}=\gamma_{\pm 1} (x^{\pm 1})^{-1} \c{L}^{\pm}{}^1_1
&\quad\mbox{for $N$ even,} \\[6pt]
\mu_a=\gamma_a r_{|a|}^{-1}r_{|a|-1}^{-1} x^{-a}
&\quad\mbox{otherwise,} 
\end{array}                                             \label{defmu}
\ee
and $\gamma_a \in \b{C}$ are normalization constants fulfilling
the conditions
\be
\begin{array}{ll}
\gamma_0 = -q^{-\frac{1}{2}} h^{-1} &\quad\mbox{for $N$ odd,} \\[6pt]
\gamma_1 \gamma_{-1}=
\left\{\begin{array}{l}
-q^{-1} h^{-2}\\
k^{-2}
\end{array}\right.
&\quad\!\begin{array}{l}
\mbox{for $N$ odd,} \\
\mbox{for $N$ even,}
\end{array}\\[8pt]
\gamma_a \gamma_{-a} =
-q^{-1} k^{-2} \omega_a \omega_{a-1} &\quad\mbox{for $a>1$}. \nonumber
\end{array}                                               \label{gamma}
\ee
$h,k,\omega_a$ are defined  as in Sections
\ref{qspaces}, \ref{heisenberg}.
On the other hand, the images of $\varphi^+$ on the positive 
FRT generators read
\be
\varphi^+(\c{L}^+{}^i_j)=g^{ih}[\bar\mu_h,x^k]_{q^{-1}}g_{kj},
                                                     \label{imagel+}
\ee
where
\be
\begin{array}{ll}
\bar\mu_0=\bar\gamma_0 (x^0)^{-1}
&\quad\mbox{for $N$ odd,} \\[6pt]
\bar\mu_{\pm 1} = 
\bar\gamma_{\pm 1} (x^{\pm 1})^{-1} \c{L}^{\mp}{}^1_1
&\quad\mbox{for 
$N$ even,} \\[6pt]
\bar\mu_a = 
\bar\gamma_a r_{|a|}^{-1}r_{|a|-1}^{-1} x^{-a}
&\quad\mbox{otherwise,} 
\end{array}                                         \label{defbarmu}
\ee
and $\bar\gamma_a \in \b{C}$ normalization constants fulfilling
the conditions
\be
\begin{array}{ll}
\bar \gamma_0 = q^{\frac{1}{2}} h^{-1}
&\quad\mbox{for $N$ odd,} \\[6pt]
\bar \gamma_1 \bar \gamma_{-1} =
\left\{
\begin{array}{l}
-q h^{-2}\\
k^{-2}
\end{array}
\right.
&\quad\!\begin{array}{l}
\mbox{for $N$ odd,} \\
\mbox{for $N$ even,}
\end{array}\\[8pt]
\bar \gamma_a \bar \gamma_{-a}=
-q k^{-2} \omega_a \omega_{a-1}
&\quad\mbox{for $a>1$}.         \nonumber
\end{array}                                             \label{bargamma}
\ee
Incidentally, for odd $N$
one can choose the free parameters $\gamma_a,\bar\gamma_a$
in such a way that $\varphi^+,\varphi^-$ can be `glued' into
an algebra homomorphism $\varphi:\b{R}_q^N\cocross\uqs\to\b{R}_q^N$
\cite{CerFioMad00}. 

\bigskip

We give the explicit expression for 
$\varphi^{\pm}(\c{L}^{\pm}{}^i_j)$ in the case $N=3$:
\be
[\varphi^-\!(\!\c{L}^-{}^i_j\!)]\!=\!\!\left[\ba{ccc}
-qh\gamma_1(x^0)^{-1} r \:& & \\
q^{\frac 12}(q\!+\!1)(x^0)^{-1}x^+ \:& 1 & \\
q^{\frac 12}(q\!+\!1)(h\gamma_1rx^0)^{-1}(x^+)^2 \:\:& (1\!+\!q^{-1}) 
(h\gamma_1 r)^{-1}\:\: &
-\!(qh\gamma_1r)^{-1}x^0\!\! 
\ea\right]
\label{vf-} 
\ee
and
\be
[\varphi^+(\c{L}^+{}^i_j)]\!=\!\left[\ba{ccc}
-h\bar\gamma_1 r^{-1}x^0\:\: & q^{-\frac 12}\bar\gamma_1 k r^{-1} x^-\: & 
q^{-2} k \bar\gamma_1(rx^0)^{-1}(x^-)^2\\
& 1 & q^{-\frac 12} (q^{-1}+1)(x^0)^{-1}x^- \\
& & -(h\bar\gamma_1x^0)^{-1}r
\ea\right]
\label{vf+} 
\ee

\bigskip

When $q\in\b{R}^+$ the real structure of $\b{R}_q^N$ is given by
\be
(x^i)^*=x^jg_{ji}.                                  \label{qreal}
\ee
Note that when $N$ is odd
$\mu_0,\bar\mu_0$, which are completely determined by their definitions,
are such that $\mu_0^*=-q^{-1}\bar\mu_0$. We fix the other 
$\gamma_a,\bar\gamma_a$ so that for any $a$
\be
\mu_a^*=-q^{-1}g_{ab}\bar\mu_b.                   \label{qreal'}
\ee
This was already considered in Ref. \cite{CerFioMad00} and
requires
\be
\ba{ll}
\gamma_{\pm 1}^*=-\bar\gamma_{\mp 1}\qquad & \qquad\mbox{if $N$ even}\\
\gamma_a^*=-\bar\gamma_{-a}\cases{\mbox{1 if }\:\:\:a<0 \cr q^{-2}
\mbox{ if }\: a>0} & \qquad\mbox{otherwise.}
\ea
\ee
As a consequence,
\bea
\left[\varphi^-(\c{L}^-{}^i_j)\right]^* &\stackrel{(\ref{imagel-})}{=}&
\left(g^{ih}[\mu_h,x^k]_qg_{kj}\right)^* \nn
&\stackrel{(\ref{qreal})}{=} & g^{ih}\,[x^j,\mu_h^*]_q\nn           
&\stackrel{(\ref{qreal'})}{=} & [\bar\mu_i,x^j]_{q^{-1}}\nn
&\stackrel{(\ref{imagel+})}{=} & g_{ih}\varphi^+(\c{L}^+{}^h_k)g^{kj}\nn
&\stackrel{(\ref{qreal0})}{=} & 
\varphi^+\left[(\c{L}^-{}^i_j)^*\right]\nonumber
\eea
In other words
\be
[\varphi^{\pm}(g)]^*=\varphi^{\mp}(g^*).
\ee
\medskip

When $|q|=1$
\be
(x^i)^*=x^i  \qquad                    \label{|q|=1'}
\ee
Note that when $N$ is odd
$\mu_0,\bar\mu_0$, which are completely determined by their definitions,
are such that $\mu_0^*=-q\mu_0=\bar\mu_0$. We fix the other 
$\gamma_a,\bar\gamma_a$ so that for any $a$
\be
\mu_a^*=-q\mu_a, \qquad\qquad\bar\mu_a^*=-q^{-1}\bar\mu_a.       
\label{|q|=1"}
\ee
This requires
\be
\ba{ll}
\gamma_{\pm 1}^*=-\gamma_{\pm 1}\qquad & \qquad\mbox{if $N$ even}\\
\gamma_a^*=-\gamma_a\cases{\mbox{1 if }\:\:\:a<0 \cr q^{-2}
\mbox{ if }\: a>0} & \qquad\mbox{otherwise.}
\ea
\ee
As a consequence,
\bea
\left[\varphi^-(\c{L}^-{}^i_j)\right]^* &\stackrel{(\ref{imagel-})}{=}&
\left(g^{ih}[\mu_h,x^k]_qg_{kj}\right)^* \nn
&\stackrel{(\ref{|q|=1'})}{=} & -q^{-1}g^{hi}\,[\mu_h^*,x^k]_q\,g_{jk}\nn           
&\stackrel{(\ref{|q|=1"})}{=} & g^{hi}\,[\mu_h,x^k]_q\,g_{jk}\nn           
&\stackrel{(\ref{defU}),(\ref{imagel-})}{=} & 
U^{-1}{}^i_r\varphi^-(\c{L}^-{}^r_s)U^i_r\nn
&\stackrel{(\ref{|q|=1})}{=} & 
\varphi^-\left[(\c{L}^-{}^i_j)^*\right].\qquad \nonumber
\eea
Similarly one proves that 
$[\varphi^-(\c{L}^-{}^i_j)]^*=\varphi^-[(\c{L}^-{}^i_j)^*]$.
In other words, $\varphi^{\pm}$ are $*$-homomorphisms.

\subsect{The maps $\varphi$ for the deformed
Heisenberg algebras}
\label{appstar}

In Ref. \cite{Fiocmp95} we constructed an algebra homomorphism
$\varphi:\uqs\cross\c{A}_1\to \c{A}_1$, where $\c{A}_1$
denotes the \uqs-covariant (deformed) Heisenberg
algebra, such that $\varphi$ is a $*$-homomorphism 
\be
\varphi(g^*)=\varphi(g)^*                  \label{*hom}
\ee
on the compact
section of \uqs (what requires $q\in\b{R}^+$).
One can easily prove the same result also for the noncompact section
(\ref{|q|=1}) of $\g=so(N)$ as well as the compact and noncompact
sections of $\g=sl(N)$. This can be done maybe most rapidly using as
a set of generators the socalled "vector fields'' $Z^i_j$ 
\cite{SchWatZum93},
which are related to the FRT generators by
\be
Z^i_j=\c{L}^+{}^i_h S\c{L}^-{}^h_j.                  \label{defZ}
\ee
{From} (\ref{qReal}), (\ref{|q|=1}) one immediately finds
\bea
&& (Z^i_j)^*= Z^j_i \qquad\qquad \qquad\qquad \mbox{if }q\in\b{R}^+  \\
&& (Z^i_j)^*= U^{-1}{}^i_a\,(S^{-1}\c{L}^-{}^h_b)\,\c{L}^+{}^a_h\, U^b_j
\qquad\qquad \mbox{if }|q|=1;
\eea
if $\g=so(N)$ the second relation reduces to
\be
(Z^i_j)^*= U^{-1}{}^a_bZ^b_c\hat R^{-1}{}^{ci}_{aj}.
\ee
In Ref. \cite{ChuZum95} the explicit expression of
$\varphi(Z^i_j)$ in terms of the $x$'s and $\partial$'s
is given both for $g=sl(N)$ and $\g=so(N)$, and it is not difficult
to show that on these generators (and therefore on all of
\uqg) (\ref{*hom}) is satisfied. In performing the calculations
one has to keep in mind that the authors of Ref. \cite{ChuZum95}
work with the left action, rather than with the right, so
one has to switch to the conventions described in section
\ref{leftsec}, but, as explained there, this wil not modify
the result (\ref{*hom}). As an intermediate step, we give the 
action of the $*$-structure on the coordinates and derivatives
for the case $\g=so(N)$, in the notation used there:
\bea
&& (x_h)^*=g^{hk}x_k, \:\: (\partial^i)^*=-q^{-N}\hat\partial_i 
\qquad\qquad \mbox{if }q\in\b{R}^+\\
&& (x_h)^*=x_h, \:\: (\partial^i)^*=-q^N U^{-1}{}^i_j\partial^j,
\:\: (\hat\partial_i)^*=-q^{-N}\partial_i \qquad\:\mbox{if }|q|=1.
\qquad \:
\eea

\end{document}